\documentclass[a4paper,12pt]{article}

\usepackage{latexsym}
\usepackage{amssymb}
\usepackage{theorem}
\usepackage{amsmath}
\usepackage{amscd}
\usepackage[dvips]{graphicx}
\usepackage{xcolor}

\newtheorem{theorem}{Theorem}[section]
\newtheorem{corollary}[theorem]{Corollary}
\newtheorem{lemma}[theorem]{Lemma}
\newtheorem{example}[theorem]{Example}
\newtheorem{proposition}[theorem]{Proposition}
\newtheorem{remark}[theorem]{Remark}
\newtheorem{definition}[theorem]{Definition}

\newcommand{\demo}{\par\noindent{\it Proof. \/}\ }
\newcommand{\enD}{\hfill $\Box$\vspace{3truemm} \par}
\newcommand{\R}{\mathbb{R}}

\newcommand{\rank}{\operatorname{rank}}

\newcommand{\bn}{\mbox{\boldmath $n$}}
\newcommand{\bt}{\mbox{\boldmath $t$}}
\newcommand{\bs}{\mbox{\boldmath $s$}}
\newcommand{\bb}{\mbox{\boldmath $b$}}
\newcommand{\ba}{\mbox{\boldmath $a$}}

\newcommand{\bx}{\mbox{\boldmath $x$}}

\newcommand{\be}{\mbox{\boldmath $e$}}
\newcommand{\bmu}{\mbox{\boldmath $\mu$}}

\newcommand{\wgamma}{\widetilde{\gamma}}
\newcommand{\wnu}{\widetilde{\nu}}
\newcommand{\wbmu}{\widetilde{\bmu}}
\newcommand{\well}{\widetilde{\ell}}
\newcommand{\wm}{\widetilde{m}}
\newcommand{\wn}{\widetilde{n}}
\newcommand{\walpha}{\widetilde{\alpha}}

\newcommand{\wkappa}{\widetilde{\kappa}}
\newcommand{\wtau}{\widetilde{\tau}}
\newcommand{\wbb}{\widetilde{\bb}}
\newcommand{\wbt}{\widetilde{\bt}}
\newcommand{\wbn}{\widetilde{\bn}}

\newcommand{\wmu}{\widetilde{\mu}}
\newcommand{\wvarphi}{\widetilde{\varphi}}

\begin{document}

\title{Singularities of translation surfaces under the linearly dependent condition in Euclidean $3$-space}
\author{T. Fukunaga, N. Nakatsuyama and M. Takahashi}

\date{\today}

\maketitle

\begin{abstract}
We investigated singular points of translation surfaces under the linearly independent condition. 
In this paper, as completion, we investigate singular points of translation surfaces under the linearly dependent condition, using the theories of generalised framed surfaces and framed surfaces. 
We introduce the notion of translation generalised framed surfaces and investigate types of singular points that appear on translation generalised framed base surfaces. 
\end{abstract}

\renewcommand{\thefootnote}{\fnsymbol{footnote}}
\footnote[0]{2020 Mathematics Subject classification: 58K05, 57R45, 53A05}
\footnote[0]{Key Words and Phrases. translation surface, generalised framed surface, framed surface, singularity.}

\section{Introduction}

The origin of translation surfaces is found in the text of G. Darboux \cite{Darboux}.  
The surfaces later known as {\it Darboux surfaces}, are defined by the form $\Phi(s,t) = A(t) \cdot \gamma(s) + \wgamma(t)$, where $\gamma(s)$ and $\wgamma(t)$ are space curves in $\mathbb{R}^3$ and $A(t)$ is an orthogonal matrix.
The surfaces present the movement of a curve by a uniparametric family of rigid motions of the Euclidean $3$-space. 
If we set $A(t)$ to be the identity matrix, we obtain a translation surface.

We consider framed curves as space curves with singular points containing non-degenerate curves (cf. \cite{Honda-Takahashi-2016}). 
Moreover, we consider generalised framed surfaces as surfaces with singular points containing not only in framed surfaces, but also in one-parameter families of framed curves (cf. \cite{Takahashi-Yu}). 

In \cite{Fukunaga-Takahashi-2022}, we investigated singular points of the translation surfaces of framed curves under the linearly independent condition. 
As completion, we investigate singular points of translation surfaces of framed curves under the linearly dependent condition, using the theories of generalised framed surfaces and framed surfaces. 

In order to investigate translation surfaces, we introduce a frame matrix for two framed curves. 
The frame matrix with a condition is determined the frame of framed curves (Theorem \ref{property-translation-matrix}). 
Using the translation matrix, we detect the basic invariants of the translation surface as a generalised framed surface. 
We introduce the notion of translation generalised framed surfaces and investigate types of singular points that appear on translation generalised framed base surfaces. 
We give conditions for a translation surface to have cross cap and $S^{\pm}_1$ singular points (Theorems \ref{S0} and \ref{S1}). 
In the case of the linearly independent conditions, translation surfaces never have cross cap and $S^{\pm}_1$ singular points. 
Moreover, it is natural for translation surfaces with linearly dependent conditions to appear as self translation surfaces generated by one framed curve. 
Furthermore, we give a condition that the translation surface to be a framed base surface under the linearly dependent condition. 
We also investigate singular points of translation surfaces as framed base surfaces (Theorems \ref{criteria-case-I}, \ref{criteria-case-II}, \ref{criteria-case-III} and \ref{criteria-case-IV}).
\par
All maps and manifolds considered in this paper are differentiable of class $C^{\infty}$ unless stated otherwise.

\bigskip
\noindent
{\bf Acknowledgements}. 
The second author was supported by JST SPRING Grant Number JPMJSP2153.
The third author was partially supported by JSPS KAKENHI Grant Number JP 24K06728.

\section{Translation surfaces as generalised framed surfaces}

Let $\R^3$ be the $3$-dimensional Euclidean space equipped with the inner product $\ba \cdot \bb = a_1 b_1 + a_2 b_2 + a_3 b_3$, 
where $\ba = (a_1, a_2, a_3)$ and $\bb = (b_1, b_2, b_3) \in \R^3$. 
The norm of $\ba$ is given by $\vert \ba \vert = \sqrt{\ba \cdot \ba}$ and the vector product is given by 
$$
\ba \times \bb={\rm det}
\left(
\begin{array}{ccc}
\be_1 & \be_2 & \be_3\\
a_1 & a_2 & a_3\\
b_1 & b_2 & b_3
\end{array}
\right),
$$
where $\{\be_1, \be_2, \be_3\}$ is the canonical basis of $\R^3$. 
Let $S^2$ be the unit sphere in $\R^3$, that is, $S^2=\{\ba \in \R^3| |\ba|=1\}$.
We denote the $3$-dimensional smooth manifold $\{(\ba,\bb) \in S^2
\times S^2| \ba \cdot \bb=0\}$ by $\Delta$. 
Let $I$ and $J$ be intervals in $\R$.
\par
Let $(\gamma,\nu_1,\nu_2): I \rightarrow \mathbb{R}^3 \time \Delta$ and $(\widetilde\gamma,\widetilde\nu_1, \widetilde\nu_2) : J \rightarrow \mathbb{R}^3 \time \Delta$ be framed curves with curvatures $(\ell, m, n, \alpha)$ and $(\widetilde\ell, \widetilde{m}, \widetilde{n}, \widetilde{\alpha})$, respectively. 
For framed curves, see Appendix \ref{framed-curves} and \cite{Honda-Takahashi-2016}. 
It follows that 
\begin{align*}
\begin{pmatrix}
\nu_{1u}(u) \\
\nu_{2u}(u) \\
\mu_u(u) 
\end{pmatrix}
&=
\mathcal{F}(u)
\begin{pmatrix}
\nu_1(u) \\
\nu_2(u) \\
\mu(u) 
\end{pmatrix}
=
\begin{pmatrix}
0 & \ell(u) & m(u) \\
-\ell(u) & 0 & n(u) \\
-m(u) & -n(u) & 0 
\end{pmatrix}
\begin{pmatrix}
\nu_1(u) \\
\nu_2(u) \\
\mu(u) 
\end{pmatrix},\\
\begin{pmatrix}
\wnu_{1v}(v) \\
\wnu_{2v}(v) \\
\wmu_v(v) 
\end{pmatrix}
&=
\widetilde{\mathcal{F}}(v)
\begin{pmatrix}
\wnu_1(v) \\
\wnu_2(v) \\
\wmu(v) 
\end{pmatrix}
=
\begin{pmatrix}
0 & \well(v) & \wm(v) \\
-\well(v) & 0 & \wn(v) \\
-\wm(v) & -\wn(v) & 0 
\end{pmatrix}
\begin{pmatrix}
\wnu_1(v) \\
\wnu_2(v) \\
\wmu(v) 
\end{pmatrix},
\end{align*}
where $\mu=\nu_1 \times \nu_2$ and $\wmu=\wnu_1 \times \wnu_2$.

\begin{definition}\label{translation-surface}
We say that $\bx:I \times J \to \R^3$ is a {\it translation surface} if $\bx(u,v)=\gamma(u)+\wgamma(v)$.
\end{definition}

Since $\{\nu_1,\nu_2,\mu\}$ and $\{\wnu_1,\wnu_2,\wmu\}$ are moving frames of $\gamma$ and $\wgamma$, respectively, there exists a smooth mapping $T: I \times J \to SO(3)$ such that 
\begin{align}\label{frame_relation}
\begin{pmatrix}
\wnu_1(v) \\
\wnu_2(v) \\
\wmu(v) 
\end{pmatrix}
=
T(u,v)
\begin{pmatrix}
\nu_1(u) \\
\nu_2(u) \\
\mu(u) 
\end{pmatrix},
\end{align}
where 
$$
T(u,v)=(t_{ij}(u,v))_{1 \le i, j \le 3}=
\begin{pmatrix}
\wnu_1(v) \cdot \nu_1(u) & \wnu_1(v) \cdot \nu_2(u) & \wnu_1(v) \cdot \mu(u) \\
\wnu_2(v) \cdot \nu_1(u) & \wnu_2(v) \cdot \nu_2(u) & \wnu_2(v) \cdot \mu(u) \\
\wmu(v) \cdot \nu_1(u) & \wmu(v) \cdot \nu_2(u) & \wmu(v) \cdot \mu(u) 
\end{pmatrix}.
$$
We say that $T$ is a {\it frame matrix} of framed curves $(\gamma,\nu_1,\nu_2)$ and $(\wgamma,\wnu_1,\wnu_2)$. 
Here 
$SO(3)$ is the set of $3 \times 3$ special orthogonal matrices, that is,  
$SO(3)=\{A \in GL(3,\R) \ |\ ^tA A=I_3, {\rm det}A=1 \}$, where $GL(3,\R)$ is the set of $3 \times 3$ general linear matrices, $^tA$ is the transpose of the matrix $A$ and $I_3$ is the $3 \times 3$ unit matrix. 
\par
By differentiating \eqref{frame_relation}, we have the following.
\begin{proposition}\label{relation_integrable}
Under the above notations, we have 
\begin{align}\label{relation_translation}
T_u(u,v) = -T(u,v)\mathcal{F}(u), \
T_v(u,v) = \widetilde{\mathcal{F}}(v)T(u,v).
\end{align}
It follows that we have 
\begin{align*}
t_{i1u}(u,v) &= \ell(u) t_{i2}(u,v) + m(u) t_{i3}(u,v), \\
t_{i2u}(u,v) &= -\ell(u) t_{i1}(u,v) + n(u) t_{i3}(u,v), \\
t_{i3u}(u,v) &= -m(u) t_{i1}(u,v) - n(u)  t_{i2}(u,v), \\
t_{1jv}(u,v) &= \well(v) t_{2j}(u,v) + \wm(v) t_{3j}(u,v), \\
t_{2jv}(u,v) &= -\well(v) t_{1j}(u,v) + \wn(v) t_{3j}(u,v), \\
t_{3jv}(u,v) &= -\wm(v) t_{1j}(u,v) - \wn(v)  t_{2j}(u,v),
\end{align*}
for $i,j=1,2,3$.
\end{proposition}

\begin{theorem}\label{property-translation-matrix}
$(1)$ Let $(\gamma, \nu_1, \nu_2) : I \rightarrow \mathbb{R}^3 \times \Delta$ and $(\wgamma, \wnu_1, \wnu_2) : I \rightarrow \mathbb{R}^3 \times \Delta$ be framed curves. 
Then the frame matrix $T : I \times J \rightarrow SO(3)$ satisfies $T_{uv} = T_v {}^tT T_u$. 
\par
$(2)$ Suppose that $T: I \times J \rightarrow SO(3)$ is a smooth mapping with $T_{uv} = T_v {}^tT T_u $, and $\alpha : I \rightarrow \mathbb{R}^3$ and $\walpha : J \rightarrow \mathbb{R}^3$ are smooth functions.  
Then there exist framed curves $(\gamma, \nu_1, \nu_2) : I \rightarrow \mathbb{R}^3 \times \Delta$ and $(\wgamma, \wnu_1, \wnu_2) : I \rightarrow \mathbb{R}^3 \times \Delta$ such that 
$\gamma_u(u) = \alpha(u)\mu(u)$, $\wgamma_v(v) = \walpha(v)\wmu(v)$ and
the frame matrix of the framed curves $(\gamma, \nu_1, \nu_2)$ and $(\wgamma, \wnu_1, \wnu_2)$ is given by $T$. 
\end{theorem}
\demo
$(1)$ By differentiating \eqref{relation_translation}, we have 
$T_{uv}=-T_v \mathcal{F}=-\widetilde{\mathcal{F}}T \mathcal{F}$. 
Since $ ^tT T=T ^t T=I_3$, we have $\mathcal{F}=- ^tTT_u$ and $\widetilde{\mathcal{F}}={T_v} ^tT$. 
It follows that $T_{uv}=-({T_v} ^tT)T(- ^tTT_u)={T_v} ^tT T_u$.
\par
$(2)$ Suppose that $T: I \times J \rightarrow SO(3)$ is a smooth mapping with $T_{uv} = T_v {}^tT T_u $. 
We define $\mathcal{F}$ and $\widetilde{\mathcal{F}}:I \times J \to M(3,3)$ by 
$$
\mathcal{F}(u,v)=- ^tT(u,v)T_u(u,v), \ \widetilde{\mathcal{F}}(u,v)=T_v(u,v) ^tT(u,v),
$$
where $M(3,3)$ is the set of $3 \times 3$ matrices. 
By using $( ^tT)_uT+ \!^tTT_u=0$ and $( ^tT)_vT+ \!^tTT_v=0$, we have 
\begin{align*}
\mathcal{F}_v(u,v) &=-( ^tT)_v(u,v)T_u(u,v)-\! ^tT(u,v)T_{uv}(u,v)\\
&= \! ^tT(u,v)(T_v(u,v) ^tT(u,v)T_u(u,v)-T_{uv}(u,v))=0,\\
\widetilde{\mathcal{F}}_u(u,v) &=T_{vu}(u,v) ^tT(u,v)+T_{v}(u,v) ( ^tT(u,v))_u\\
&= (T_{vu}(u,v) ^t-T_v(u,v) ^tT(u,v)T_u(u,v))T(u,v)=0. 
\end{align*}
Therefore $\mathcal{F}(u,v)=\mathcal{F}(u)$ and $\widetilde{\mathcal{F}}(u,v)=\widetilde{\mathcal{F}}(v)$. 
Since  
\begin{align*}
^t\mathcal{F} &=- {} ^t(^tT T_u)=- {} ^t(T_u) T=-(^tT)_uT= ^tTT_u {} ^tT T=  {} ^tTT_u=-\mathcal{F}, \\
^t\widetilde{\mathcal{F}}&= {} ^t({T_v} {} ^tT)=T {} ^t(T_v)=T( {} ^tT)_v=-T {} ^tTT_v {} ^tT=-T_v {} ^tT=-\widetilde{\mathcal{F}},
\end{align*}
$\mathcal{F}(u)$ and $\widetilde{\mathcal{F}}(v)$ are alternative matrices.
By existence of a solution of a system of linear differential equations, there exist smooth mappings 
$(\nu_1,\nu_2,\mu): I \to SO(3)$ and $(\wnu_1,\wnu_2,\wmu):J \to SO(3)$ such that 
\begin{align*}
\begin{pmatrix}
\nu_{1u} (u)\\
\nu_{2u} (u)\\
\mu_u (u)
\end{pmatrix}
=
\mathcal{F}(u)
\begin{pmatrix}
\nu_1(u) \\
\nu_2(u) \\
\mu (u)
\end{pmatrix}, 
\begin{pmatrix}
\wnu_{1v} (v)\\
\wnu_{2v} (v)\\
\wmu_v(v)
\end{pmatrix}
=
\widetilde{\mathcal{F}}(v)
\begin{pmatrix}
\wnu_1(v) \\
\wnu_2(v) \\
\wmu (v)
\end{pmatrix},
\begin{pmatrix}
\wnu_{1}(v_0) \\
\wnu_{2}(v_0) \\
\wmu(v_0)
\end{pmatrix}
=T(u_0,v_0)
\begin{pmatrix}
\nu_1(u_0) \\
\nu_2(u_0) \\
\mu (u_0)
\end{pmatrix}
\end{align*}
for a point $(u_0,v_0) \in I \times J$.
If $\gamma:I \to \R^3$ and $\wgamma:J \to \R^3$ are given by
$$
\gamma(u)=\int \alpha(u) \mu(u) du, \ \wgamma(u)=\int \alpha(v) \wmu(v) dv,
$$
then $(\gamma, \nu_1, \nu_2) : I \rightarrow \mathbb{R}^3 \times \Delta$ and $(\wgamma, \wnu_1, \wnu_2) : I \rightarrow \mathbb{R}^3 \times \Delta$ are framed curves with $\gamma_u(u) = \alpha(u)\mu(u)$ and $\wgamma_v(v) = \walpha(v)\wmu(v)$. 
Moreover, since 
\begin{align*}
&\frac{\partial}{\partial u}
\left(
\left(\begin{array}{cccc}
\wnu_{1}\\
\wnu_{2}\\
\wmu
\end{array}\right) -
T
\left(\begin{array}{cccc}
\nu_{1}\\
\nu_{2}\\
\mu
\end{array}\right)\right)
= -T_u
\left(\begin{array}{cccc}
\nu_{1}\\
\nu_{2}\\
\mu
\end{array}\right)
- T{} \mathcal{F}
\left(\begin{array}{cccc}
\nu_{1}\\
\nu_{2}\\
\mu
\end{array}\right)
=0,\\
&
\frac{\partial}{\partial v}
\left( 
^tT
\left(\begin{array}{cccc}
\wnu_{1}\\
\wnu_{2}\\
\wmu
\end{array}\right) -
\left(\begin{array}{cccc}
\nu_{1}\\
\nu_{2}\\
\mu
\end{array}\right)\right)
= ({} ^tT)_v
\left(\begin{array}{cccc}
\wnu_{1}\\
\wnu_{2}\\
\wmu
\end{array}\right)
+ {} ^tT \widetilde{\mathcal{F}}
\left(\begin{array}{cccc}
\wnu_{1}\\
\wnu_{2}\\
\wmu
\end{array}\right)
=0,
\end{align*}
there exist smooth function $A:J \to \R^3$ and $B:I \to \R^3$ such that
$$
\left(\begin{array}{cccc}
\wnu_{1}(v)\\
\wnu_{2}(v)\\
\wmu(v)
\end{array}\right) -
T(u,v)
\left(\begin{array}{cccc}
\nu_{1}(u)\\
\nu_{2}(u)\\
\mu(u)
\end{array}\right)
=A(v), \ {}
^tT(u,v)
\left(\begin{array}{cccc}
\wnu_{1}(v)\\
\wnu_{2}(v)\\
\wmu(v)
\end{array}\right) -
\left(\begin{array}{cccc}
\nu_{1}(u)\\
\nu_{2}(u)\\
\mu(u)
\end{array}\right)
= B(u).
$$
It follows that $ ^tT(u,v) A(v)=B(u)$ for all $(u,v) \in I \times J$.
By the initial condition, we have $A(v_0)=0$ and $B(u_0)=0$.
By $ ^tT(u_0,v)A(v)=B(u_0)=0$, we have $A(v)=0$ for all $v \in J$ and hence $B(u)=0$ for all $u \in I$. 
Therefore, $T$ is the frame matrix of the framed curves $(\gamma, \nu_1, \nu_2)$ and $(\wgamma, \wnu_1, \wnu_2)$.
\enD

The translation surface $\bx$ is a generalised framed base surface, in fact, we have the following. 
For generalised framed surface, see Appendix \ref{GFS-FS} and \cite{Takahashi-Yu}.
\begin{proposition}\label{TGFS} 
Let $\bx: I \times J \to \R^3$ be the translation surface of framed curves $(\gamma,\nu_1,\nu_2)$ and $(\wgamma,\wnu_1,\wnu_2)$.
\par
$(1)$ $(\bx,\nu_1,\nu_2):I \times J \to \R^3 \times \Delta$ is a generalised framed surface with basic invariants
\begin{align*}
\begin{pmatrix}
a_1(u,v) & b_1(u,v) & c_1(u,v) \\
a_2(u,v) & b_2(u,v) & c_2(u,v) 
\end{pmatrix} &= 
\begin{pmatrix}
0 & 0 & \alpha(u) \\
\walpha(v) t_{31}(u,v) & \walpha(v) t_{32}(u,v) & \walpha(v) t_{33}(u,v)
\end{pmatrix},
\\
\begin{pmatrix}
e_1(u,v) & f_1(u,v) & g_1(u,v) \\
e_2(u,v) & f_2(u,v) & g_2 (u,v)
\end{pmatrix} &= 
\begin{pmatrix}
\ell(u) & m(u) & n(u) \\
0 & 0 & 0
\end{pmatrix},
\\
A(u,v) = -\alpha(u) \walpha(v) t_{32}(u,v),& \ B(u,v) = \alpha(u) \walpha(v) t_{31}(u,v),
\end{align*}
where $(\nu_1,\nu_2): I \times J \to \Delta$ is defined by $\nu_i(u,v)=\nu_i(u), i=1,2$.
\par
$(2)$ $(\bx, \wnu_1, \wnu_2):I \times J \to \R^3 \times \Delta$ is a generalised framed surface with basic invariants
\begin{align*}
\begin{pmatrix}
a_1(u,v) & b_1(u,v) & c_1(u,v) \\
a_2(u,v) & b_2(u,v) & c_2(u,v) 
\end{pmatrix} &= 
\begin{pmatrix}
\alpha(u) t_{13}(u,v) & \alpha(u) t_{23}(u,v) & \alpha(u) t_{33}(u,v) \\
0 & 0 & \walpha \\
\end{pmatrix},
\\
\begin{pmatrix}
e_1(u,v) & f_1(u,v) & g_1(u,v) \\
e_2(u,v) & f_2(u,v) & g_2(u,v) 
\end{pmatrix} &=
\begin{pmatrix}
0 & 0 & 0 \\
\well(v) & \wm(v) & \wn(v) 
\end{pmatrix},
\\
A(u,v) = \alpha(u) \walpha(v) t_{23}(u,v),& \ B(u,v) = -\alpha(u) \walpha(v) t_{13}(u,v),
\end{align*}
where $(\wnu_1,\wnu_2): I \times J \to \Delta$ is defined by $\wnu_i(u,v)=\wnu_i(v), i=1,2$.
\end{proposition}

\demo
$(1)$ By a direct calculation, we have
\begin{align*}
a_1(u,v) &= \bx_u(u,v) \cdot \nu_1(u,v) = \alpha(u) \mu(u,v) \cdot \nu_1(u,v) = 0, \\ 
b_1(u,v) &= \bx_u(u,v) \cdot \nu_2(u,v) = \alpha(u) \mu(u,v) \cdot \nu_2(u,v) = 0, \\
c_1(u,v) &= \bx_u(u,v) \cdot \mu(u,v) = \alpha(u) \mu(u,v) \cdot \mu(u,v) = \alpha(u), \\  
a_2(u,v) &= \bx_v(u,v) \cdot \nu_1(u,v) = \walpha(v) \wmu(u,v) \cdot \nu_1(u,v) = \walpha(v) t_{31}(u,v), \\
b_2(u,v) &= \bx_v(u,v) \cdot \nu_2(u,v) = \walpha(v) \wmu(u,v) \cdot \nu_2(u,v) = \walpha(v) t_{32}(u,v), \\
c_2(u,v) &= \bx_v(u,v) \cdot \mu(u,v) = \walpha(v) \wmu(u,v) \cdot \mu(u,v) = \walpha(v) t_{33}(u,v), \\
e_1(u,v) &= \nu_{1u}(u,v) \cdot \nu_2(u,v) = \ell(u), \ e_2(u,v) = \nu_{1v}(u,v) \cdot \nu_2(u,v) = 0, \\ 
f_1(u,v) &= \nu_{1u}(u,v) \cdot \mu(u,v) = m(u), \ \ f_2(u,v) = \nu_{1v}(u,v) \cdot \mu(u,v) = 0, \\ 
g_1(u,v) &= \nu_{2u}(u,v) \cdot \mu(u,v) = n(u), \ \ g_2(u,v) = \nu_{2v}(u,v) \cdot \mu(u,v) = 0. 
\end{align*}
Moreover, we have
\begin{align*}
A(u,v) &= \det \begin{pmatrix} 0 & \walpha(v) t_{32}(u,v) \\ \alpha(u) & \walpha(v) t_{33}(u,v) \end{pmatrix} = -\alpha(u) \walpha(v) t_{32}(u,v), \\ 
B(u,v) &= - \det \begin{pmatrix} 0 & \alpha(u) t_{31}(u,v)\\ \alpha(u) & \walpha(v) t_{33}(u,v) \end{pmatrix} = \alpha(u) \walpha(v) t_{31}(u,v).
\end{align*}
\par
$(2)$ By a similar calculation, we have the result.
\enD
By Proposition \ref{TGFS}, 
$p=(u,v)$ is a singular point of $\bx$ if and only if the one of conditions satisfies 
\begin{align*}
{\rm(i)}  \ \alpha(u)=0, \quad 
{\rm(ii)}  \ \walpha(v)=0, \quad 
{\rm(iii)}  \ t_{31}(p)=t_{32}(p)=0. 
\end{align*}
The cases $\rm(i)$ and $\rm(ii)$ are corresponding to the singular point of $\gamma$ and $\wgamma$, respectively. 
Moreover, $\mu(u) \times \wmu(v)=0$ if and only if the condition $\rm(iii)$ $t_{31}(p)=t_{32}(p)=0$. 
On the other hand, if $\mu(u) \times \wmu(v) \not=0$ at $p$, the we can define unit normal vector field of $\bx$ around $p$. 
Hence $\bx$ is a frontal (or a framed base surface) around $p$. 
We say that the {\it linearly dependent condition} (respectively, {\it linearly independent condition}) of framed curves $(\gamma,\nu_1,\nu_2)$ and $(\wgamma,\wnu_1,\wnu_2)$ at $p$ if $\mu(u) \times \wmu(v)=0$ (respectively, $\mu(u) \times \wmu(v) \not=0$). 
We have investigated singular points of translation surfaces under the linearly independent condition in \cite{Fukunaga-Takahashi-2022}. 

By Proposition \ref{TGFS} and Theorem \ref{FBS.condition}, we have the following. 

\begin{proposition}\label{FBS.condition-translation}
Let $(\bx, \nu_1, \nu_2) : I \times J \rightarrow \mathbb{R}^3 \times \Delta$ be a translation generalised framed surface. 
\par
$(1)$ If $\bx$ is a framed base surface, then $t_{31}$ and $t_{32}$ (respectively, $t_{13}$ and $t_{23}$) are linearly dependent. 
\par
$(2)$ Suppose that the set of regular points of $\bx$ is dense in $I \times J$. 
If $t_{31}$ and $t_{32}$ (respectively, $t_{13}$ and $t_{23}$) are linearly dependent, then $\bx$ is a framed base surface.
\end{proposition}


\subsection{A criteria for the cross cap ($S_0$ singularity)}\label{S0-singularities}
The symbols and settings in this subsection continue to be those used above. 
If $\alpha(u_0)=0$ or $\walpha(v_0)=0$, that is, $u_0$ is a singular point of $\gamma$ or $v_0$ is a singular point of $\wgamma$, then $p_0=(u_0,v_0)$ is not a cross cap singular point of translation surface $\bx$ by using criteria of cross cap singular point (\cite{Whitney1943}). See also \cite{Fukunaga-Takahashi-2022}. 
\par
Consider a singular point $p_0 \in I \times J$ of a generalised translation framed base surface $\bx : I \times J \rightarrow \mathbb{R}^3$.
Suppose that $\alpha(u_0)\walpha(v_0) \not=0$. 
In this case, the vector field $\eta = -\walpha(\partial/\partial u) + \alpha t_{33}(\partial / \partial v)$ satisfies $d\bx_{p_0}(\eta_{p_0}) = 0$, that is, $\ker(d\bx_{p_0}) = \langle \eta_{p_0} \rangle_{\mathbb{R}}$. Here, we use the relation $\mu(u_0) = t_{33}(p_0)\wmu(v_0)$ derived from $\alpha(u_0)\walpha(v_0) \neq 0$ and $\rank d\bx_{p_0}< 2 $. If $\xi$ is defined as $\partial/\partial u$, then $\xi_{p_0}$ and $\eta_{p_0}$ are linearly independent since $\alpha(u_0)\walpha(v_0) \neq 0$.

We define the function $\varphi : I \times J \rightarrow \mathbb{R}$ by $\varphi(u,v) = \det(\xi\bx,\eta\bx,\eta\eta\bx)(u,v)$. By a direct calculation, we have
\begin{align*}
\varphi &= -\alpha^3\walpha^3 t_{33}\left(-m t_{32}+n t_{31}\right)-\alpha^4\walpha^2 t_{33}^3\left( \wm t_{23}-\wn t_{13} \right),\\
\varphi_u &= -3\alpha^2 \alpha_u \walpha^3t_{33}\left(-m t_{32}+n t_{31}\right) 
- \alpha^3 \walpha^3 t_{33u}\left(-m t_{32}+n t_{31}\right)\\
&\quad -\alpha^3 \walpha^3 t_{33} \left( -m_u t_{32} - m t_{32u} +n_u t_{31} +n t_{31u} \right)\\
&\quad -4\alpha^3 \alpha_u \walpha^2 t_{33}^3\left( \wm t_{23}-\wn t_{13} \right) - 3\alpha^4\walpha^2t_{33}^2t_{33u}(\wm t_{23}- \wn t_{13})\\
&\quad -\alpha^4\walpha^2 t_{33}^3 \left\{ \wm (-m t_{21} - n t_{22}) - \wn(-m t_{11} - n t_{12})\right\}, \\
\varphi_v &= -3\alpha^2\walpha^2\walpha_v t_{33}(-m t_{32}+n t_{31})
-\alpha^3 \walpha^3 t_{33v}\left( -m t_{32} + n t_{31} \right)\\ 
&\quad -\alpha^3 \walpha^3 t_{33}\left( -m t_{32v} + n t_{31v} \right)\\
&\quad -2\alpha^4 \walpha \walpha_v t_{33}^3\left( \wm t_{23} - \wn t_{13} \right) - 3\alpha^4\walpha^2t_{33}^2t_{33v}(\wm t_{23}- \wn t_{13})\\
&\quad -\alpha^4\walpha^2 t_{33}^3 \left(\wm_v t_{23} +  \wm t_{23v} -\wn_vt_{13} -\wn t_{13v}\right).
\end{align*}
\begin{theorem}\label{S0}
Let $(\gamma, \nu_1, \nu_2) : I \rightarrow \mathbb{R}^3 \times \Delta$ and $(\wgamma, \wnu_1, \wnu_2) : J \rightarrow \mathbb{R}^3 \times \Delta$ be framed curves with curvatures $(\ell, m, n, \alpha)$ and $(\well, \wm, \wn, \walpha)$, respectively.
$\bx(u,v)=\gamma(u)+\wgamma(v)$ at $p_0=(u_0,v_0)$ is a cross cap (an $S_0$ singular point) if and only if $\alpha(u_0)\walpha(v_0) \neq 0$, $|t_{33}(p_0)| = 1$ and 
\begin{eqnarray*}
m(u_0) (\wm(v_0) t_{21}(p_0) - \wn(v_0) t_{11}(p_0)) + n(u_0) (\wm(v_0) t_{22}(p_0) - \wn(v_0) t_{12}(p_0)) \neq 0. 
\end{eqnarray*}
\end{theorem}
\demo
Using Theorem \ref{CMM}, calculate $\xi\varphi$. Since $\xi = \partial / \partial u$, 
\begin{align*}
\xi\varphi_{p_0} &= \varphi_u(p_0) \\ &= 
-3\alpha(u_0)^2 \alpha_u(u_0) \walpha(v_0)^3t_{33}(p_0)\left(-m(u_0) t_{32}(p_0)+n(u_0) t_{31}(p_0)\right) 
\\&\quad - \alpha(u_0)^3 \walpha(v_0)^3 t_{33u}(p_0)\left(-m(u_0) t_{32}(p_0)+n(u_0) t_{31}(p_0)\right)\\
&\quad -\alpha(u_0)^3 \walpha(v_0)^3 t_{33}(p_0) \left( -m_u(u_0) t_{32}(p_0) - m(u_0) t_{32u}(p_0) +n_u(u_0) t_{31}(p_0) +n(u_0) t_{31u}(p_0) \right)\\
&\quad -4\alpha(u_0)^3 \alpha_u(u_0) \walpha(v_0)^2 t_{33}(p_0)^3\left( \wm(v_0) t_{23}(p_0)-\wn(v_0) t_{13}(p_0) \right) \\
&\quad- 3\alpha(u_0)^4\walpha(v_0)^2t_{33}(p_0)^2t_{33u}(p_0)(\wm(v_0) t_{23}(p_0)- \wn(v_0) t_{13}(p_0))\\
&\quad -\alpha(u_0)^4\walpha(v_0)^2 t_{33}(p_0)^3 \left\{ \wm(v_0) (-m(u_0) t_{21}(p_0) - 
n(u_0) t_{22}(p_0)) \right. \\& \left. \quad - \wn(v_0)(-m(u_0) t_{11}(p_0) - n(u_0) t_{12}(p_0))\right\}.
\end{align*}
From the definition of the frame matrix and $p_0$ is a singular point, $t_{33}(p_0)^2 = 1$ and $t_{31}(p_0) = t_{32}(p_0)=t_{13}(p_0)=t_{23}(p_0)=0$. Combining this with $\alpha(u_0)\walpha(v_0)\neq 0$, we obtain $\xi\varphi_{p_0} \neq 0$ if and only if $m(u_0) (\wm(v_0) t_{21}(p_0) - \wn(v_0) t_{11}(p_0)) + n(u_0) (\wm(v_0) t_{22}(p_0) - \wn(v_0) t_{12}(p_0)) \neq 0$. Using the above and what was stated at the beginning of this subsection, we obtain this theorem.
\enD
\begin{example}[Cross cap]\label{example-S0}
{\rm
Let $(\gamma, \nu_1, \nu_2) : \R \rightarrow \mathbb{R}^3 \times \Delta$ and  
$(\wgamma, \wnu_1, \wnu_2) : \R \rightarrow \mathbb{R}^3 \times \Delta$ be 
\begin{align*}
\gamma(u) &= \left(u, \frac{u^2}{2}, 0\right), \ \nu_1(u) = \frac{1}{\sqrt{1+u^2}}\left( -u,1, 0 \right), \ \nu_2(u) =  \left( 0, 0, 1 \right), \\
\wgamma(v) &= \left(v, 0, \frac{v^2}{2} \right), \ \wnu_1(v) = \frac{1}{\sqrt{1+v^2}}\left( v, 0, -1 \right), \ \wnu_2(v) =  \left( 0 ,1, 0 \right).
\end{align*} 
Then $(\gamma,\nu_1,\nu_2)$ and $(\wgamma,\wnu_1,\wnu_2)$ are framed curves with curvatures 
\begin{align*}
(\ell(u), m(u), n(u), \alpha(u)) &= \left(0, \frac{-1}{1+u^2}, 0, \sqrt{1+u^2}  \right), \\
(\well(v), \wm(v), \wn(v), \walpha(v)) &= \left(0, \frac{1}{1+v^2}, 0, \sqrt{1+v^2}  \right). 
\end{align*}
By a direct calculation, we have 
\begin{eqnarray*}
\mu(u) = \frac{1}{\sqrt{1+u^2}}\left(1, u ,0 \right), \ \wmu(v) = \frac{1}{\sqrt{1+v^2}}\left(1, 0, v \right).
\end{eqnarray*}
It follows that $(\bx, \nu_1, \nu_2): \R \times \R \to \R^3 \times \Delta$, 
\begin{align*}
&\bx(u,v)=\gamma(u)+\wgamma(v)=\left(u+v,\frac{u^2}{2},\frac{v^2}{2}\right), \\
&\nu_1(u,v)=\nu_1(u)=\frac{1}{\sqrt{1+u^2}}\left( -u,1, 0 \right), \ \nu_2(u,v)=\nu_2(u) =\left( 0, 0, 1 \right), 
\end{align*}
is a generalised framed surface. 
By a direct calculation, the frame matrix is given by
$$
T(u,v)=\begin{pmatrix}
\frac{-uv}{\sqrt{(1+u^2)(1+v^2)}} &\frac{-1}{\sqrt{1+v^2}} & \frac{v}{\sqrt{(1+u^2)(1+v^2)}} \\
\frac{1}{\sqrt{1+u^2}} & 0 & \frac{u}{\sqrt{1+u^2}} \\
\frac{-u}{\sqrt{(1+u^2)(1+v^2)}} & \frac{v}{\sqrt{1+v^2}} & \frac{1}{\sqrt{(1+u^2)(1+v^2)}}
\end{pmatrix}.
$$
Therefore, we have $\alpha(0)\walpha(0)=1$, $t_{33}(0,0)=1$ and 
\begin{eqnarray*}
m(0) (\wm(0) t_{21}(0,0) - \wn(0) t_{11}(0,0)) + n(0)(\wm(0) t_{22}(0,0) - \wn(0) t_{12}(0,0))=-1 \neq 0. 
\end{eqnarray*}
By Theorem \ref{S0}, $(0,0)$ is a cross cap singular point of the translation surface $\bx$. 
\begin{figure}[h!]
\begin{center}
\includegraphics[width=55mm,height=55mm]{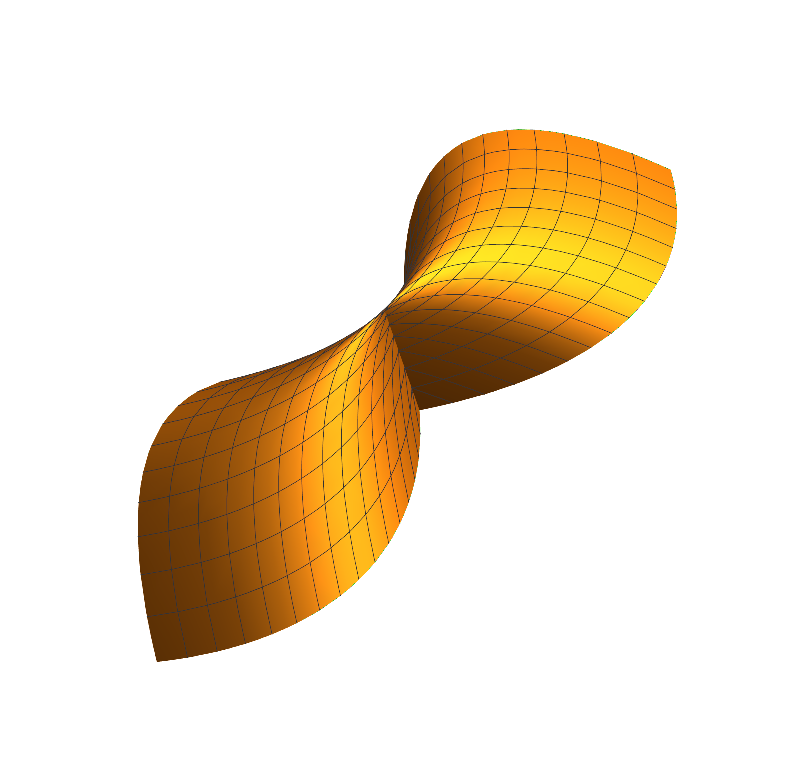}
\caption{The translation surface $\bx$ in Example \ref{example-S0}. $(0,0)$ is a cross cap singular point.}
\end{center}
\end{figure}
}
\end{example}
\begin{remark}{\rm
It is easy to see that $(0,0)$ is a cross cap singular point of $\bx$ in Example \ref{example-S0} by diffeomorphisms of the domain and of the target ($\mathcal{A}$-equivalent). 
In fact, 
\begin{eqnarray*}
\bx(u,v) = \left(u+v, \frac{u^2}{2}, \frac{v^2}{2} \right) 
 = \left(U, \frac{(U-v)^2}{2}, \frac{v^2}{2} \right) 
 \sim_{\mathcal{A}} \left( U, Uv, v^2 \right).
\end{eqnarray*}
}
\end{remark}

We consider $\gamma:I \to \R^3$ and $\wgamma:J \to \R^3$ are non-degenerate curves with arc-length parameter.
Then we have the Frenet frame $\{\bt,\bn,\bb \}$ of $\gamma$ with curvature $\kappa$ and torsion $\tau$, and the Frenet frame  $\{\wbt,\wbn,\wbb \}$ of $\wgamma$ with curvature $\wkappa$ and torsion $\wtau$, respectively.
It follows that $(\gamma,\bn,\bb):I \to \R^3 \times \Delta$ and $(\wgamma,\wbn,\wbb):I \to \R^3 \times \Delta$ are framed curves with curvatures $(\ell,m,n,\alpha)=(\tau,-\kappa,0,1)$ and 
$(\well,\wm,\wn,\walpha)=(\wtau,-\wkappa,0,1)$, respectively (cf. Appendix \ref{framed-curves}). 
\begin{corollary}\label{S0frenet}
Suppose that $\gamma:I \to \R^3$ and $\wgamma:J \to \R^3$ are non-degenerate curves with arc-length parameter. 
Then $\bx(u,v)=\gamma(u)+\wgamma(v)$ at $p_0=(u_0,v_0)$ is a cross cap (an $S_0$ singular point) if and only if $|t_{33}(p_0)|=1$ and $t_{21}(p_0) \neq 0$.
\end{corollary}
\demo
Applying Theorem \ref{S0} to case $(\ell,m,n,\alpha)=(\tau,-\kappa,0,1)$ and $(\well,\wm,\wn,\walpha)=(\wtau,-\wkappa,0,1)$, we obtain \begin{eqnarray*}
m(u_0) (\wm(v_0) t_{21}(p_0) - \wn(v_0) t_{11}(p_0)) + n(u_0) (\wm(v_0) t_{22}(p_0) - \wn(v_0) t_{12}(p_0)) = \kappa(u_0)\wkappa(v_0)t_{21}(p_0). 
\end{eqnarray*}
Thus, $\kappa(u_0)\wkappa(v_0) \neq 0$ yields the assertion.
\enD

\subsection{Criterion for the $S^\pm_1$ singularities}\label{S1-singularities}
Let $p_0 = (u_0, v_0)$ be a singular point of $\bx$. 
If $\alpha(u_0)=0$ or $\walpha(v_0)=0$, that is, $u_0$ is a singular point of $\gamma$ or $v_0$ is a singular point of $\wgamma$, then $(u_0,v_0)$ is not a $S_{1}^{\pm}$ singular point of translation surface $\bx$ by using criteria of $S^\pm_1$ singular point (\cite{Whitney1943}). See also \cite{Fukunaga-Takahashi-2022}.
\par
Suppose that $\alpha(u_0)\walpha(v_0) \not =0$ and $\mu(u_0) = t_{33}(p_0)\wmu(v_0)$. 
We use the same notations $\xi, \eta$ and $\varphi$ as in section \ref{S0-singularities}. 
By using the relational expressions $t_{13}(p_0)=t_{23}(p_0) = t_{31}(p_0)=t_{32}(p_0)=0$, $t_{33}(p_0)^2 = 1$, $t_{31u}(p_0) = m(u_0)t_{33}(p_0)$ and
$t_{31uu}(p_0) = \left(\ell(u_0) n(u_0) +m_u(u_0) \right) t_{33}(p_0)$, we obtain
\begin{align*}
\varphi_{uu}(p_0) &= \alpha(u_0)^3\walpha(v_0)^2\left(-\walpha(v_0) (\ell(u_0) (m(u_0)^2 + n(u_0)^2) - n(u_0) m_u(u_0)  + m(u_0)n_u(u_0)) \right. \\ 
&\quad + t_{33}(p_0) \alpha(u_0) (\ell(u_0) (n(u_0) (-t_{21}(p_0) \wm(v_0) + t_{11}(p_0)\wn(v_0))\\ 
&\quad + m(u_0) (t_{22}(p_0) \wm(v_0) - t_{12}(p_0) \wn(v_0))) + (t_{21}(p_0) \wm(v_0) - t_{11}(p_0) \wn(v_0))m_u(u_0) \\ 
&\quad + (t_{22}(p_0) \wm(v_0) - t_{12}(p_0) \wn(v_0)) n_u(u_0)) + 8 t_{33}(p_0) ((m(u_0) t_{21}(p_0) + n(u_0) t_{22}(p_0)) \wm(v_0) \\ 
&\quad \left.- (m(u_0) t_{11}(p_0) + n(u_0)t_{12}(p_0)) \wn(v_0))\alpha_u(u_0)\right),\\
\varphi_{uv}(p_0) &= t_{33}(p_0) \walpha(v_0) \alpha(u_0)^2 (\walpha(v_0)^2 
\alpha(u_0) (\ell(u_0) (m(u_0) (t_{11}(p_0) \wm(v_0) + t_{21}(p_0) \wn(v_0)) \\
&\quad +n(u_0) (t_{12}(p_0) \wm(v_0) + t_{22}(p_0) \wn(v_0))) - (t_{12}(p_0) \wm(v_0) + t_{22}(p_0) \wn(v_0)) m_u(u_0) \\
&\quad + (t_{11}(p_0) \wm(v_0) + t_{21}(p_0) \wn(v_0)) n_u(u_0)) \\
&\quad - \walpha(v_0) \alpha(u_0)^2 (m(u_0) (t_{21}(p_0) (\well(v_0) \wn(v_0) - \wm_v(v_0)) \\
&\quad + t_{11}(p_0) (\well(v_0) \wm(v_0) + \wn_v(v_0))) + n(u_0) (t_{22}(p_0) (\well(v_0) \wn(v_0) - \wm_v(v_0)) \\
&\quad + t_{12}(p_0) (\well(v_0) \wm(v_0) + \wn_v(v_0)))) + 2 ((m(u_0) t_{21}(p_0) + n(u_0) t_{22}(p_0)) \wm(v_0) \\
&\quad - (m(u_0) t_{11}(p_0) + n(u_0) t_{12}(p_0)) \wn(v_0)) \alpha(u_0)^2 \walpha_v(v_0)\\
&\quad +  3 (n(u_0) (t_{11}(p_0) \wm(v_0) + t_{21}(p_0) \wn(v_0)) - m(u_0) (t_{12}(p_0) \wm(v_0) \\
&\quad +t_{22}(p_0) \wn(v_0))) \walpha(v_0)^2 \alpha_u(u_0)),\\
\varphi_{vv}(p_0) &=
t_{33}(p_0) \walpha(v_0)^2 \alpha(u_0)^3 (t_{33}(p_0) \alpha(u_0) (\well(v_0) (\wm(v_0)^2 + \wn(v_0)^2) - \wn(v_0) \wm_v(v_0) \\
&\quad + \wm(v_0) \wn_v(v_0)) + n(u_0) (\walpha(v_0) (t_{11}(p_0) (-\well(v_0) \wn(v_0) + \wm_v(v_0)) \\
&\quad + t_{21}(p_0) (\well(v_0) \wm(v_0) + \wn_v(v_0))) + 6 (t_{11}(p_0) \wm(v_0) + t_{21}(p_0) \wn(v_0)) \walpha_v(v_0)) \\
&\quad + m(u_0) (-\walpha(v_0) (t_{12}(p_0) (-\well(v_0) \wn(v_0) + \wm_v(v_0))  \\
&\quad + t_{22}(p_0) (\well(v_0) \wm(v_0) + \wn_v(v_0))) - 6 (t_{12}(p_0) \wm(v_0) + t_{22}(p_0) \wn(v_0)) \walpha_v(v_0))).
\end{align*}

\begin{theorem}\label{S1}
Let $(\gamma, \nu_1, \nu_2) : I \rightarrow \mathbb{R}^3 \times \Delta$ and $(\wgamma, \wnu_1, \wnu_2) : J \rightarrow \mathbb{R}^3 \times \Delta$ be framed curves with curvatures $(\ell, m, n, \alpha)$ and $(\well, \wm, \wn, \walpha)$, respectively. 
\par
$(1)$ $\bx(u,v)=\gamma(u)+\wgamma(v)$ at $p_0=(u_0,v_0)$ is an $S_1^+$ singular point if and only if $\alpha(u_0)\walpha(v_0) \neq 0$, $|t_{33}(p_0)| = 1$, 
$m(u_0) (\wm(v_0) t_{21}(p_0) - \wn(v_0) t_{11}(p_0)) + n(u_0) (\wm(v_0) t_{22}(p_0) - \wn(v_0) t_{12}(p_0)) = 0,$
$\varphi_{uu}(p_0)\varphi_{vv}(p_0) - \varphi_{uv}(p_0)^2 < 0$ and 
\begin{align*}
&(\walpha(v_0)m(u_0)+\alpha(u_0)(\wm(v_0)t_{11}(p_0)+\wn(v_0)t_{21}(p_0)), \\ 
& \ \walpha(v_0)n(u_0)+\alpha(u_0)(\wm(v_0)t_{12}(p_0)+\wn(v_0)t_{22}(p_0)) \not=(0,0).
\end{align*}
\par
$(2)$ $\bx(u,v)=\gamma(u)+\wgamma(v)$ at $p_0=(u_0,v_0)$ is an $S_1^-$ singular point if and only if $\alpha(u_0)\walpha(v_0) \neq 0$, $|t_{33}(p_0)| = 1$, 
$m(u_0) (\wm(v_0) t_{21}(p_0) - \wn(v_0) t_{11}(p_0)) + n(u_0) (\wm(v_0) t_{22}(p_0) - \wn(v_0) t_{12}(p_0)) = 0$ and 
$\varphi_{uu}(p_0)\varphi_{vv}(p_0) - \varphi_{uv}(p_0)^2 > 0$.
\end{theorem}
\demo
$(1)$ \  By a direct calculation,
\begin{align*}
\xi\bx(p_0) &= \alpha(u_0)\mu(u_0), \\
\eta\eta\bx(p_0) &= \walpha(v_0)^2\left\{ \alpha_u(u_0)\mu(u_0) + 
\alpha(u_0)(-m(u_0) \nu_1(u_0) -n(u_0) \nu_2(u_0) ) -\alpha_u(u_0)t_{33}(p_0)\wmu(v_0) \right\} 
\\
&\quad +
\alpha(u_0)^2\left\{ -\walpha_v(v_0)t_{33}(p_0)\mu(u_0) + \walpha_v(v_0)\wmu(v_0)+\walpha(v_0)(-\wm(v_0)\wnu_1(v_0) - \wn(v_0)\wnu_2(v_0)) \right\}.
\end{align*}
Since $\mu(u_0)$ and $\wmu(v_0)$ are linearly independent,  $\xi\bx(p_0)$ and $\eta\eta\bx(p_0)$ are linearly independent if and only if 
\begin{align*}
&\bigl(\walpha(v_0)m(u_0)+\alpha(u_0)(\wm(v_0)t_{11}(p_0)+\wn(v_0)t_{21}(p_0)), \\ 
& \ \walpha(v_0)n(u_0)+\alpha(u_0)(\wm(v_0)t_{12}(p_0)+\wn(v_0)t_{22}(p_0)\bigr) \not=(0,0).
\end{align*}
From the definition of the frame matrix and $\alpha(u_0)\walpha(v_0) \neq 0$, $p_0$ is a singular point of $\bx$ if and only if $t_{33}(p_0)^2 = 1$. 
By the assumption $m(u_0) (\wm(v_0) t_{21}(p_0) - \wn(v_0) t_{11}(p_0)) + n(u_0) (\wm(v_0) t_{22}(p_0) - \wn(v_0) t_{12}(p_0)) = 0$, $\varphi$ has a critical point at $p_0$. 
From the above discussion and Theorem \ref{CMM} $(2)$, we obtain the assertion. 

$(2)$ \  By the above discussion and Theorem \ref{CMM} $(3)$, we obtain the assertion.

\enD

\begin{example}[$S^+_1$ singularity]\label{example-S1+}{\rm
Let $(\gamma, \nu_1, \nu_2) : \R \rightarrow \mathbb{R}^3 \times \Delta$ and  
$(\wgamma, \wnu_1, \wnu_2) : \R \rightarrow \mathbb{R}^3 \times \Delta$ be
\begin{align*}
\gamma(u) &= \left(u, \frac{u^3}{3}, 0\right), \ \nu_1(u) = \frac{1}{\sqrt{1+u^4}}\left(-u^2, 1, 0 \right), \ \nu_2(u) =  \left( 0, 0, 1 \right), \\
\wgamma(v) &= \left(v, 0, \frac{v^2}{2} \right), \ \wnu_1(v) = \frac{1}{\sqrt{1+v^2}}\left( v, 0, -1 \right), \ \wnu_2(v) =  \left(0 ,1, 0 \right).
\end{align*} 
Then $(\gamma,\nu_1,\nu_2)$ and $(\wgamma, \wnu_1, \wnu_2)$ are framed curves with curvatures 
\begin{align*}
(\ell(u), m(u), n(u), \alpha(u)) &= \left(0, \frac{-2u}{1+u^4}, 0, \sqrt{1+u^4}  \right), \\
(\well(v), \wm(v), \wn(v), \walpha(v)) &= \left(0, \frac{1}{1+v^2}, 0, \sqrt{1+v^2}  \right).
\end{align*}
By a direct calculation, we have 
$$
\mu(u) = \frac{1}{\sqrt{1+u^4}}\left(1, u^2 ,0 \right), \ \wmu(v) = \frac{1}{\sqrt{1+v^2}}\left(1, 0, v \right).
$$
It follows that $(\bx,\nu_1,\nu_2):\R \times \R \to \R^3 \times \Delta$, 
\begin{align*}
&\bx(u,v) = \gamma(u) + \wgamma(v)=\left(u+v,\frac{u^3}{3},\frac{v^2}{2}\right), \\
&\nu_1(u,v)=\nu_1(u)=\frac{1}{\sqrt{1+u^4}}\left(-u^2, 1, 0 \right), \ \nu_2(u,v)=\nu_2(u) =\left( 0, 0, 1 \right) 
\end{align*}
is a generalised framed surface.
By a direct calculation, the frame matrix is given by
$$
T(u,v)=\begin{pmatrix}
\frac{-u^2v}{\sqrt{(1+u^4)(1+v^2)}} &\frac{-1}{\sqrt{1+v^2}} & \frac{v}{\sqrt{(1+u^4)(1+v^2)}} \\
\frac{1}{\sqrt{1+u^4}} & 0 & \frac{u^2}{\sqrt{1+u^4}} \\
\frac{-u^2}{\sqrt{(1+u^4)(1+v^2)}} & \frac{v}{\sqrt{1+v^2}} & \frac{1}{\sqrt{(1+u^4)(1+v^2)}}
\end{pmatrix}.
$$
Moreover, we have $t_{12}(0,0) = -1$, 
$\displaystyle m_u(u) = \frac{6u^4 -2}{(1+u^4)^2} $, $\displaystyle\alpha_u(u) = \frac{2u^3}{\sqrt{1+u^4}}$, 
$\displaystyle\wm_v(v) = \frac{-2v}{(1+v^2)^2}$, $\displaystyle\walpha_v(v) = \frac{v}{\sqrt{1+v^2}}$. 
By a direct calculation, $\alpha(0)\walpha(0) \not=0$, $t_{33}(0,0)=1$, $m(0) (\wm(0) t_{21}(0,0) - \wn(0) t_{11}(0,0)) + n(0) (\wm(0) t_{22}(0,0) - \wn(0) t_{12}(0,0)) = 0$, $\det (\text{Hess}\varphi(0,0)) = -4 < 0$ and 
\begin{align*}
&(\walpha(0)m(0)+\alpha(0)(\wm(0)t_{11}(0,0)+\wn(0)t_{21}(0,0)), \\ 
& \ \walpha(0)n(0)+\alpha(0)(\wm(0)t_{12}(0,0)+\wn(0)t_{22}(0,0))=(0,1).
\end{align*} 
By Theorem \ref{S1}, $(0,0)$ is an $S_1^+$ singular point of the translation surface $\bx$.

\begin{figure}[h!]
\begin{center}
\includegraphics[height=55mm]{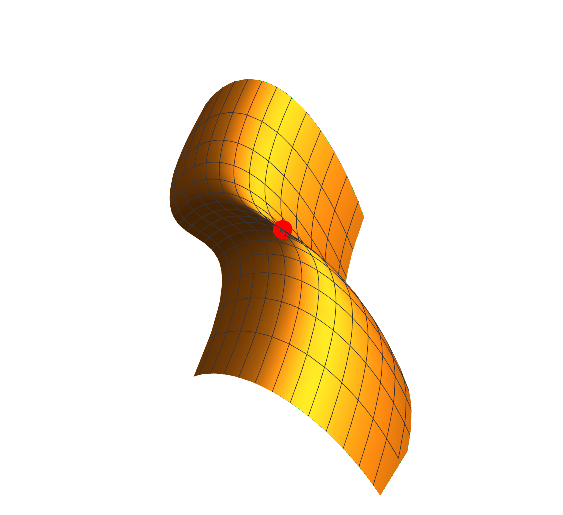}
\caption{The translation surface $\bx$ in Example \ref{example-S1+}. $(0,0)$ is an $S_1^+$ singular point and the dot in the figure is the point $\bx(0,0)$.}
\end{center}
\end{figure}
}
\end{example}
We give an example of $S^-_1$ singularity in Example \ref{example-S1-}. 

\begin{corollary}\label{S1frenet}
Let $\gamma : I \rightarrow \mathbb{R}^3$ and $\widetilde{\gamma}:J \to \R^3$ be non-degenerate curves with arc-length parameter. 
Suppose that $|t_{33}(p_0)|=1$. 
\par
$(1)$ $p_0$ is an $S_1^+$ singular point of $\bx$ if and only if 
$t_{21}(p_0) = 0$, $\kappa(u_0)+\wkappa(v_0)t_{11}(p_0) \neq 0$ and 
$$
\tau(u_0)\wtau(v_0)(\kappa(u_0)^2+\wkappa(v_0)^2)t_{11}(p_0)
-\kappa(u_0)\wkappa(v_0)(\tau(u_0)^2+\wtau(v_0)^2)<0.
$$
\par
$(2)$ $p_0$ is an $S_1^-$ singular point of $\bx$ if and only if $t_{21}(p_0) = 0$ and 
$$
\tau(u_0)\wtau(v_0)(\kappa(u_0)^2+\wkappa(v_0)^2)t_{11}(p_0)
-\kappa(u_0)\wkappa(v_0)(\tau(u_0)^2+\wtau(v_0)^2)>0.
$$
\end{corollary}
\demo
$(1)$ \ By substituting the relations $(\ell, m, n, \alpha) = (\tau, -\kappa, 0, 1)$ and $(\well, \wm, \wn, \walpha) = (\wtau, -\wkappa, 0, 1)$ into Theorem \ref{S1}, we obtain
\begin{align*}
\varphi_{uu}(p_0) &= \kappa_u(u_0)\wkappa(v_0)t_{21}(p_0)t_{33}(p_0) + \kappa(u_0)\tau(u_0)\wkappa(v_0)t_{22}(p_0)t_{33}(p_0) - \tau(u_0)\kappa(u_0)^2, \\
\varphi_{uv}(p_0) &= \left\{ \kappa(u_0)\tau(u_0)\wkappa(v_0)-\kappa(u_0)\wkappa(v_0)\wtau(v_0))t_{11}(p_0) \right. \\ 
&\quad \left.
- \kappa_u(u_0)\wkappa(v_0)t_{12}(p_0)+ \kappa(u_0)\wkappa_v(v_0)t_{21}(p_0)\right\}t_{33}(p_0), \\
\varphi_{vv}(p_0) &= -\kappa(u_0)\wkappa_v(v_0)t_{12}(p_0)t_{33}(p_0) - \kappa(u_0)\wkappa(v_0)\wtau(v_0)t_{22}(p_0)t_{33}(p_0) + \wkappa(v_0)^2 \wtau(v_0), \\
\eta\eta\bx(p_0) &= (\kappa(u_0) + \wkappa(v_0) t_{11}(p_0))\nu_1(u_0)+ \wkappa(v_0) t_{12}(p_0) \nu_2(u_0),\\
\xi\bx(p_0) &= \mu(u_0).
\end{align*}
By the assumption $t_{21}(p_0)=0$, the equality $t_{12}(p_0)=0$ follows. Thus, $\xi\bx(p_0)$ and $\eta\eta\bx(p_0)$ are linearly independent if and only if $\kappa(u_0)+\wkappa(v_0)t_{11}(p_0) \neq 0$.
Then, from $|t_{33}(p_0)| = 1$ and $T(p_0) \in SO(3)$, we have
\begin{eqnarray*}
\begin{pmatrix} 
t_{11}(p_0) & t_{12}(p_0)\\
t_{21}(p_0) & t_{22}(p_0)
\end{pmatrix} 
\in O(2)
\end{eqnarray*}
and 
\begin{align*}
&\det(\text{Hess}\varphi(p_0)) = 
\det\begin{pmatrix} 
\varphi_{uu}(p_0) & \varphi_{uv}(p_0) \\
\varphi_{vu}(p_0) & \varphi_{vv}(p_0)
\end{pmatrix} \\
&= \det\begin{pmatrix} 
\kappa_u\wkappa t_{21}t_{33} + \kappa\tau\wkappa t_{22}t_{33}-\kappa^2\tau  & 
\{(\kappa \tau \wkappa - \kappa \wkappa \wtau)t_{11}-\kappa_u\wkappa t_{12}+\kappa\wkappa_v t_{21}\}t_{33} \\
\{(\kappa \tau \wkappa -\kappa \wkappa \wtau)t_{11}-\kappa_u\wkappa t_{12}+\kappa\wkappa_v t_{21}\}t_{33}  & -\kappa\wkappa_v t_{12}t_{33} - \kappa\wkappa\wtau t_{22}t_{33}+\wkappa^2\wtau
\end{pmatrix} (p_0)\\
&= \kappa(u_0)\tau(u_0)\wkappa(v_0)\wtau(v_0)\left( \kappa(u_0)^2 +\wkappa(v_0)^2 \right)t_{11}(p_0) - \kappa(u_0)^2\wkappa(v_0)^2\left(\tau(u_0)^2 + \wtau(v_0)^2 \right).
\end{align*}
Thus, $\det(\text{Hess}\varphi(p_0)) < 0$ if and only if
\begin{eqnarray*}
\tau(u_0)\wtau(v_0)\left( \kappa(u_0)^2 + \wkappa(v_0)^2 \right)t_{11}(p_0)
-\kappa(u_0)\wkappa(v_0)\left( \tau(u_0)^2  + \wtau(v_0)^2 \right) < 0.
\end{eqnarray*}
By above and Theorem \ref{S1} (1), we obtain the assertion.

$(2)$ \ $\det (\text{Hess}\varphi(p_0)) > 0$ if and only if
\begin{eqnarray*}
\tau(u_0)\wtau(v_0)\left( \kappa(u_0)^2 + \wkappa(v_0)^2 \right)t_{11}(p_0)
-\kappa(u_0)\wkappa(v_0)\left( \tau(u_0)^2  + \wtau(v_0)^2 \right) > 0.
\end{eqnarray*}
By above and Theorem \ref{S1} (2), we obtain the assertion.
\enD

\begin{example}[$S^-_1$ singularity]\label{example-S1-}{\rm
Let $\gamma$ and $\wgamma : \mathbb{R} \rightarrow \mathbb{R}^3$,
$$
\gamma(u) = \begin{pmatrix} 
\frac{3\sqrt{10}}{10} & 0 & -\frac{\sqrt10}{10}\\
0 & 1 & 0 \\
\frac{\sqrt{10}}{10} & 0 & \frac{3\sqrt10}{10} \\
\end{pmatrix} 
\begin{pmatrix} 
\frac{u}{\sqrt{2}} \\
\frac{\cos \sqrt{2}u}{2} \\
\frac{\sin \sqrt{2}u}{2} 
\end{pmatrix}, \ 
\wgamma(v) = 
\begin{pmatrix}
\frac{v}{\sqrt{5}} \\
\frac{2\cos \sqrt{5}v}{5} \\
\frac{2\sin \sqrt{5}v}{5} 
\end{pmatrix}
$$
be non-degenerate curves with the curvatures and torsions  $\kappa(u) = 1, \tau(u) = 1$, $\wkappa(v) = 2, \wtau(v) = 1$, respectively. 
By a direct calculation, we have 
\begin{align*}
\bt(u)&=\left(\frac{3-\cos\sqrt{2}u}{\sqrt{20}},-\frac{\sin\sqrt{2}u}{\sqrt{2}},\frac{1+3\cos\sqrt{2}u}{\sqrt{20}}\right), \\
\bn(u)&=\left(\frac{\sin\sqrt{2}u}{\sqrt{10}},-\cos\sqrt{2}u,-\frac{3\sin\sqrt{2}u}{\sqrt{10}}\right),\\
\bb(u)&=\left(\frac{3+\cos\sqrt{2}u}{\sqrt{20}},\frac{\sin\sqrt{2}u}{\sqrt{2}},\frac{1-3\cos\sqrt{2}u}{\sqrt{20}}\right),
\end{align*}
\begin{align*}
\wbt(v)&=\left(\frac{1}{\sqrt{5}},-\frac{2\sin\sqrt{5}v}{\sqrt{5}},\frac{2\cos\sqrt{5}v}{\sqrt{5}}\right), \\
\wbn(v)&=\left(0,-\cos\sqrt{5}v,-\sin\sqrt{5}v\right), \\
\wbb(v)&=\left(\frac{2}{\sqrt{5}},\frac{\sin\sqrt{5}v}{\sqrt{5}},-\frac{\cos\sqrt{5}v}{\sqrt{5}}\right).
\end{align*}
Then the frame matrix $T(u,v)=(t_{ij}(u,v))$ is given by
\begin{align*}
t_{11}(u,v)&=\wbn(v)\cdot\bn(u)=\cos\sqrt{2}u\cos\sqrt{5}v+\frac{3\sin\sqrt{2}u\sin\sqrt{5}v}{\sqrt{10}},\\
t_{12}(u,v)&=\wbn(v)\cdot\bb(u)=-\frac{\sin\sqrt{2}u\cos\sqrt{5}v}{\sqrt{2}}-\frac{\sin\sqrt{5}v(1-3\cos\sqrt{2}u)}{\sqrt{20}},\\
t_{13}(u,v)&=\wbn(v)\cdot\bt(u)=\frac{\sin\sqrt{2}u\cos\sqrt{5}v}{\sqrt{2}}-\frac{\sin\sqrt{5}v(1+3\cos\sqrt{2}u)}{\sqrt{20}},\\
t_{21}(u,v)&=\wbb(v)\cdot\bn(u)=\frac{\sqrt{2}\sin\sqrt{2}u}{5}-\frac{\cos\sqrt{2}u\sin\sqrt{5}v}{\sqrt{10}}+\frac{3\sin\sqrt{2}u\cos\sqrt{5}v}{5\sqrt{2}},\\
t_{22}(u,v)&=\wbb(v)\cdot\bb(u)=\frac{3+\cos\sqrt{2}u}{5}+\frac{\sin\sqrt{2}u\sin\sqrt{5}v}{\sqrt{10}}-\frac{\cos\sqrt{5}v(1-3\cos\sqrt{2}u)}{10},\\
t_{23}(u,v)&=\wbb(v)\cdot\bt(u)=\frac{3-\cos\sqrt{2}u}{5}-\frac{\sin\sqrt{2}u\sin\sqrt{5}v}{\sqrt{10}}-\frac{\cos\sqrt{5}v(1+3\cos\sqrt{2}u)}{10},\\
t_{31}(u,v)&=\wbt(v)\cdot\bn(u)=\frac{\sin\sqrt{2}u}{5\sqrt{2}}+\frac{2\cos\sqrt{2}u\sin\sqrt{5}v}{\sqrt{5}}-\frac{3\sqrt{2}\sin\sqrt{2}u\cos\sqrt{5}v}{5},\\
t_{32}(u,v)&=\wbt(v)\cdot\bb(u)=\frac{3+\cos\sqrt{2}u}{10}-\frac{2\sin\sqrt{2}u\sin\sqrt{5}v}{\sqrt{10}}+\frac{\cos\sqrt{5}v(1-3\cos\sqrt{2}u)}{5},\\
t_{33}(u,v)&=\wbt(v)\cdot\bt(u)=\frac{3-\cos\sqrt{2}u}{10}+\frac{2\sin\sqrt{2}u\sin\sqrt{5}v}{\sqrt{10}}+\frac{\cos\sqrt{5}v(1+3\cos\sqrt{2}u)}{5}.
\end{align*}
Therefore, we have $t_{33}(0,0) = 1$,  $t_{21}(0,0) = 0$ and
$$
\tau(0)\wtau(0)\left(\kappa(0)^2 + \wkappa(0)^2 \right) t_{11}(0,0) -\kappa(0)\wkappa(0)\left( \tau(0)^2 + \wtau(0)^2 \right) 
= 1.
$$
By Corollary \ref{S1frenet}, $(0,0)$ is an $S_1^-$ singular point of the translation surface $\bx$.
}
\end{example}
\begin{figure}[h!]
\begin{center}
\includegraphics[width=55mm,height=55mm]{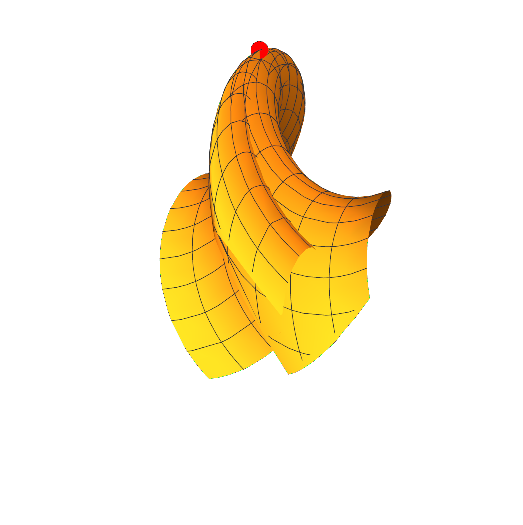}
\caption{The translation surface $\bx$ in Example \ref{example-S1-}. $(0,0)$ is an $S_1^-$ singular point and the dot in the figure is the point $\bx (0,0)$.}
\end{center}
\end{figure}

\section{Self translation surfaces}

Let $(\gamma, \nu_1, \nu_2) : I \rightarrow \mathbb{R}^3 \times \Delta$ be a framed curve with curvature $(\ell,m,n,\alpha)$. 
We consider $(\bx^{\pm}, \nu_1, \nu_2) : I \times I \rightarrow \mathbb{R}^3 \times \Delta$ as
$$
\bx^{\pm}(u,v) = \frac{\gamma(u) \pm \gamma(v)}{2},\ \nu_1(u,v)=\nu_1(u),\ \nu_2(u,v)=\nu(u).
$$
We say that $\bx^{\pm}$ is a {\it self translation surface} (cf. \cite{Akamine}). 
Since $(-\gamma, \nu_1, \nu_2):I \to \R^3 \times \Delta$ is also a framed curve with curvature $(\ell,m,n,-\alpha)$, $(\bx^{\pm}, \nu_1, \nu_2) : I \times I \rightarrow \mathbb{R}^3 \times \Delta$ is a generalised framed surface with basic invariants
\begin{align*}
\begin{pmatrix}
a_1(u,v) & b_1(u,v) & c_1(u,v) \\
a_2(u,v) & b_2(u,v) & c_2(u,v)
\end{pmatrix}
&=
\begin{pmatrix}
0 & 0 & \pm\frac{\alpha(u)}{2}  \\
\pm \frac{\alpha(v)t_{31}(u,v)}{2}  & 
\pm \frac{\alpha(v)t_{32}(u,v)}{2} &  
\pm \frac{\alpha(v)t_{33}(u,v)}{2}
\end{pmatrix},
\\
\begin{pmatrix}
e_1(u,v) & f_1(u,v) & g_1(u,v) \\
e_2(u,v) & f_2(u,v) & g_2(u,v)
\end{pmatrix}
&=
\begin{pmatrix}
\ell(u) & m(u) & n(u)  \\
0 & 0 & 0
\end{pmatrix}.
\end{align*}
Note that the frame matrix satisfies $T(u,u) = I_3$ for all $u \in I$ and the tangency condition of $(\gamma,\nu_1,\nu_2)$ and $(\pm \gamma, \nu_1,\nu_2)$ always satisfies at $(u,u) \in I \times I$. 

By Theorem \ref{S0}, we have the following corollary.

\begin{corollary}\label{self-translation-S0}
Let $(\gamma,\nu_1,\nu_2) : I \rightarrow \mathbb{R}^3 \times \Delta$ be a framed curve and $p_0 = (u_0,v_0) \in I \times I$. 
Suppose that $\mu(u_0) \times \mu(v_0)=0$. 
Then $p_0 \in I \times I$ is a cross cap singular point of $\bx^{+}$ if and only if $p_0 \in I \times I$ is a cross cap singular point of $\bx^{-}$.
\end{corollary}
\demo
Considering the application of Theorem \ref{S0} to $\bx^{+}$ and $\bx^-$, from the choice of frame, $(\well, \wm, \wn)$ in Theorem \ref{S0} is equal to $(\ell, m, n)$ in both cases $\bx^+$ and $\bx^-$. Furthermore, a frame matrix $T$ is equal in both cases $\bx^+$ and $\bx^-$. Since only $(\ell, m, n)$, $(\well, \wm, \wn)$ and $T$ appear in criterion for $S_0$ singular point, 
$p_0 \in I \times I$ is a cross cap singular point of $\bx^{+}$ if and only if $p_0 \in I \times I$ is a cross cap singular point of $\bx^{-}$.
\enD

\begin{corollary}\label{self-translation-S0-Frenet}
Let $\gamma : I \rightarrow \mathbb{R}^3$ be a non-degenerate curve with arc-length parameter and $p_0 = (u_0,v_0) \in I \times I$. 
Suppose that $\bt(u_0) \times \bt(v_0)=0$.
\par
$(1)$ $p_0 \in I \times I$ is a cross cap singular point of $\bx^\pm$ if and only if $t_{21}(p_0) \not=0$.
\par
$(2)$ Suppose that $u_0=v_0$, that is, $p_0=(u_0,u_0)$. 
Then $p_0$ is neither a cross cap singular point of $\bx^+$ nor of $\bx^-$.
\end{corollary}
\demo
$(1)$ \ Substituting the relations $(\ell, m, n, \alpha) = (\tau, -\kappa, 0, 1)$ and $(\well, \wm, \wn, \walpha) = (\tau, -\kappa, 0, \pm 1)$ into Corollary \ref{S0} completes the proof.

$(2)$ \ In both cases $\bx^+$ and $\bx^-$, $t_{21}(u_0,u_0) = \nu_2(u_0) \cdot \nu_1(u_0) = 0$. By $(1)$, $p_0$ is neither a cross cap singular point.
\enD

\begin{example}\label{sin}
{\rm 
Let $(\gamma,\nu_1,\nu_2): [-\pi, \pi] \rightarrow \mathbb{R}^3 \times \Delta$, 
\begin{align*}
\gamma(u) &= \left( \sin u, -\cos u, -\frac{1}{2}\cos2u \right), \\
\nu_1(u) &= (-\sin u, \cos u, 0), \\
\nu_2(u) &= \frac{1}{\sqrt{\sin^2 2u +1}}\left(-\sin 2u \cos u, -\sin 2u \sin u, 1  \right)
\end{align*}
be a framed curve with curvature 
$$
(\ell(u), m(u), n(u), \alpha(u)) = \left( \frac{\sin 2u}{\sqrt{\sin^2 2u + 1}}, \frac{-1}{\sqrt{\sin^2 2u + 1}}, \frac{-2\cos 2u}{\sqrt{\sin^2 2u + 1}}, \sqrt{\sin^2 2u +1} \right).
$$
By a direct calculation, we have
$$
\mu(u) = \frac{1}{\sqrt{\sin^2 2u +1}}\left(\cos u, \sin u, \sin 2u  \right).
$$
Then the frame matrix $T: [-\pi, \pi] \times  [-\pi, \pi] \to SO(3)$ of $(\gamma,\nu_1,\nu_2)$ and $(\gamma,\nu_1,\nu_2)$ is given by 
\begin{eqnarray*}
T(u,v) = 
\begin{pmatrix}
\cos(u-v) & \frac{\sin 2u \sin (v-u)}{\sqrt{\sin^2 2u + 1}} & \frac{\sin(u-v)}{\sqrt{\sin^2 2u + 1}} \\
 \frac{\sin 2v \sin (u-v)}{\sqrt{\sin^2 2v + 1}} &  \frac{\sin 2u \sin 2v \cos(u-v)+1}{\sqrt{(\sin^2 2u+1)(\sin^2 2v+1)}} & \frac{\sin 2u -\sin 2v \cos(u-v)}{\sqrt{(\sin^2 2u+1)(\sin^2 2v+1)}}\\
 \frac{\sin(v-u)}{\sqrt{\sin^2 2v + 1}} & 
 \frac{\sin 2v - \sin 2u \cos(v-u)}{\sqrt{(\sin^2 2u+1)(\sin^2 2v+1)}} & \frac{\cos(u-v)+\sin 2u\sin2v}{\sqrt{(\sin^2 2u+1)(\sin^2 2v+1)}}
\end{pmatrix}.
\end{eqnarray*}
The basic invariants of the self translation surfaces $(\bx^{\pm}, \nu_1, \nu_2)$ are given by 
\begin{align*}
\begin{pmatrix}
a_1 & b_1 & c_1 \\
a_2 & b_2 & c_2
\end{pmatrix}
&=
\begin{pmatrix}
0 & 0 & \pm\sqrt{\sin^2 2u +1}  \\
\pm \sin(u-v)  & \pm \frac{-\sin 2u \sin 2v \cos(u+v)}{\sqrt{\sin^2 2u+1}} & \pm \frac{\cos(u-v)+\sin 2u\sin2v}{\sqrt{\sin^2 2u+1}}
\end{pmatrix},\\
\begin{pmatrix}
e_1 & f_1 & g_1 \\
e_2 & f_2 & g_2
\end{pmatrix}
&=
\begin{pmatrix}
\frac{\sin 2u}{\sqrt{\sin^2 2u + 1}} & \frac{-1}{\sqrt{\sin^2 2u + 1}} & \frac{-2\cos 2u}{\sqrt{\sin^2 2u + 1}}  \\
0 & 0 & 0
\end{pmatrix}.
\end{align*}
The set of singular points $S(\bx^\pm)$ of $\bx^{\pm} : [-\pi, \pi] \times [-\pi, \pi] \rightarrow \mathbb{R}^3$ is given by 
\begin{eqnarray*}
S(\bx^{\pm}) = \left\{ (u,v) | u=v \right\} \cup \left\{ (0,\pi), \left(\frac{\pi}{2}, \frac{3\pi}{2}\right),  (\pi, 0),  \left(\frac{3\pi}{2}, \frac{\pi}{2}\right) \right\}.
\end{eqnarray*}
We put 
$$
\displaystyle \widetilde{S}(\bx^{\pm}) = \left\{ (0,\pi), \left(\frac{\pi}{2}, \frac{3\pi}{2}\right), (\pi, 0),  \left(\frac{3\pi}{2}, \frac{\pi}{2}\right) \right\}.
$$
By a direct calculation, we have $\xi\varphi = -4 \neq 0$ at $p_0 \in \widetilde{S}(\bx^{\pm})$.
By Theorem \ref{S0}, $p_0 \in \widetilde{S}(\bx^{\pm})$ is a cross cap singular point of $\bx^{\pm}$.
Since 
$$
\bx^{+}(0,\pi) = \bx^{+}(\pi,0), \ \bx^{+}\left(\frac{\pi}{2}, \frac{3\pi}{2}\right) = \bx^{+}\left(\frac{3\pi}{2}, \frac{\pi}{2}\right),
$$ 
the numbers of image of singular points are $\# \bx^{+}(\widetilde{S}(\bx_{+})) = 2$ and  $\#\bx^{-}(\widetilde{S}(\bx^{-}))=4$. 
Moreover, the set $\bx^{+}(\{(u,v)  | u=v \})$ is boundary of the image of the surface $\bx^+$ and the set $\bx^{-}(\{(u,v)  | u=v \})$ is the origin.

\begin{figure}[h!]
\begin{center}
\includegraphics[width=55mm,height=55mm]{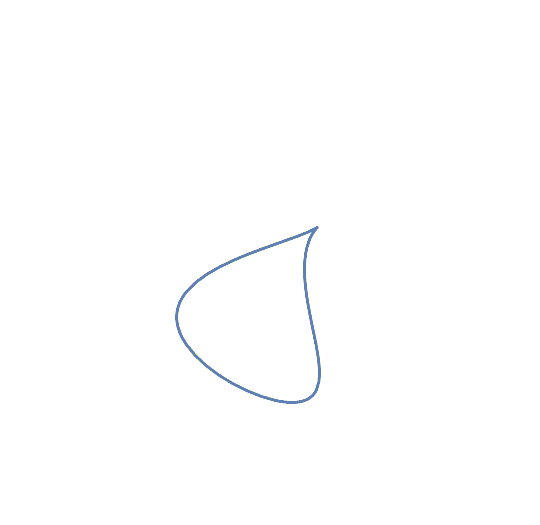}
\includegraphics[width=55mm,height=55mm]{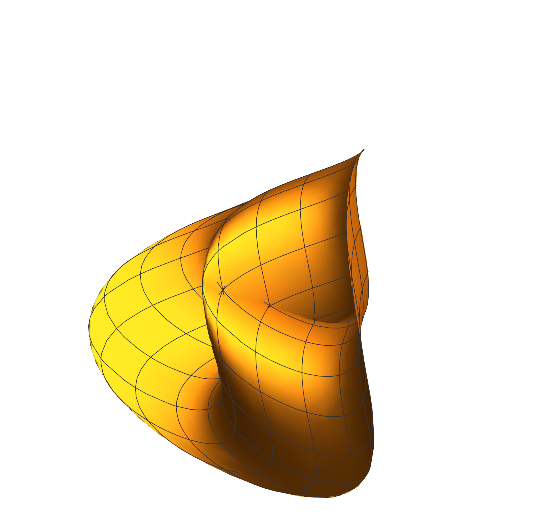}
\includegraphics[width=55mm,height=55mm]{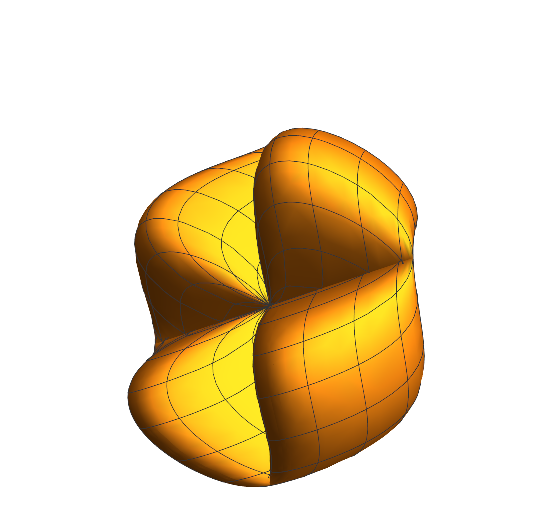}
\caption{The curve $\gamma$ and the self translation surfaces $\bx^+$ and $\bx^-$ in Example \ref{sin}.}
\end{center}
\end{figure}
}
\end{example}

In the case of $\gamma$ is a non-degenerate curve with arc-length parameter, 
we consider $S_1^{\pm}$ singularities of the self translation surface $\bx^\pm$. 
\begin{corollary}\label{self-translation-S1}
Let $\gamma : I \rightarrow \mathbb{R}^3$ be a non-degenerate curve with arc-length parameter and $p_0 = (u_0,v_0) \in I \times I$. 
Then we have the following.
\par
$(1)$ $p_0$ is an $S_1^+$ singular point of $\bx^{\pm}$ if and only if $t_{21}(p_0)=0$, 
$\kappa(u_0) \pm\kappa(v_0)t_{11}(p_0) \neq 0$ and
$$
\tau(u_0)\tau(v_0)(\kappa(u_0)^2 + \kappa(v_0)^2)t_{11}(p_0)-\kappa(u_0)\kappa(v_0)(\tau(u_0)^2+\tau(v_0)^2)<0.
$$
\par
$(2)$ $p_0$ is an $S_1^-$ singular point of $\bx^{\pm}$ if and only if $t_{21}(p_0)=0$ and 
$$
\tau(u_0)\tau(v_0)(\kappa(u_0)^2 + \kappa(v_0)^2)t_{11}(p_0)-\kappa(u_0)\kappa(v_0)(\tau(u_0)^2+\tau(v_0)^2)>0.
$$
\par
$(3)$ Suppose that $u_0=v_0$, that is, $p_0=(u_0,u_0)$. 
Then $p_0$ is neither an $S^\pm_1$ singular point of $\bx^+$ nor of $\bx^-$.
\end{corollary}
\demo
$(1)$ \ Substituting the relations $(\ell, m, n, \alpha) = (\tau, -\kappa, 0, 1)$ and $\left(\well, \wm, \wn, \walpha\right) = (\tau, -\kappa, 0, \pm 1)$ into the conditions of Theorem \ref{S1}, $\xi\bx^{\pm}(p_0)$ and $\eta\eta\bx^{\pm}(p_0)$ are linearly independent if and only if $\kappa(u_0) \pm\kappa(v_0)t_{11}(p_0) \neq 0$.
Moreover, we obtain $\det (\text{Hess} \varphi(p_0)) < 0$ if and only if
\begin{align*}
&\tau(u_0)\wtau(v_0)\left(\wkappa(v_0)t_{22}(p_0)\mp\kappa(u_0)t_{33}(p_0)\right)\left(\wkappa(v_0)t_{33}(p_0)\mp\kappa(u_0)t_{22}(p_0)\right) \\
&-\kappa(u_0)\wkappa(v_0)\left(\tau(u_0) \mp \wtau(v_0)\right)^2 < 0.
\end{align*}
In either case of the double signs, it is equivalent to
$$
\tau(u_0)\tau(v_0)(\kappa(u_0)^2 + \kappa(v_0)^2)t_{11}(p_0)-\kappa(u_0)\kappa(v_0)(\tau(u_0)^2+\tau(v_0)^2)<0.
$$

$(2)$ \ By above calculation, $\det (\text{Hess}\varphi(p_0)) > 0$ if and only if 
$$
\tau(u_0)\tau(v_0)(\kappa(u_0)^2 + \kappa(v_0)^2)t_{11}(p_0)-\kappa(u_0)\kappa(v_0)(\tau(u_0)^2+\tau(v_0)^2)>0.
$$

$(3)$ \ Substitute $p_0=(u_0,u_0)$ into $\det (\text{Hess}\varphi(p_0))$ above, $\det (\text{Hess}\varphi(p_0)) = 0$ is obtained. By Theorem \ref{S1}, the assertion follows.
\enD
\begin{corollary}\label{self-translation-S1-}
Under the same assumptions as in Corollary \ref{self-translation-S1}, suppose that $u_0 \not=v_0$ and $p_0=(u_0,v_0) \in I \times I$.
Then $p_0$ is an $S^-_1$ singular point of $\bx^{+}$ if and only if $p_0$ is an $S^-_1$ singular point of $\bx^{-}$.
\end{corollary}

\begin{remark}\label{remark-S1assumption}{\rm 
The assumption of the arc-length parameter in Corollary \ref{self-translation-S1} can be weakened as follows: $\alpha(u_0) = 1, \alpha(v_0) = \pm1$ and $\alpha_u(u_0)=\alpha_u(v_0)=0$, where $\alpha$ is the curvature of the framed curve $(\gamma, \bn, \bb)$. The examples below use this assumption. }
\end{remark}


\begin{example}\label{selfS1+exam}{\rm
$(1)$ \ Consider the regular curve $\gamma :  I \rightarrow \mathbb{R}^3$ defined by 
\begin{eqnarray*}
\gamma(u) = \left( \frac{\sin u + \sin 2u}{2\sqrt{2}} - \frac{\sin 3u}{6\sqrt{2}}, -\frac{\cos 2u}{2\sqrt{2}}, \frac{\sin 2u}{2\sqrt{2}} \right).
\end{eqnarray*}  
By a direct calculation, we have 
\begin{align*}
\gamma^{\prime}(u) &= \left( \frac{\sin u \sin 2u + \cos 2u}{\sqrt{2}}, \frac{1}{\sqrt{2}}\sin 2u, \frac{1}{\sqrt{2}}\cos 2u \right),\\
\gamma^{\prime\prime}(u) &= \left( \frac{\cos u \sin 2u + 2\sin u \cos 2u - 2\sin 2u}{\sqrt{2}}, \frac{2}{\sqrt{2}}\cos 2u, -\frac{2}{\sqrt{2}}\sin 2u \right),\\
\gamma^{\prime\prime\prime}(u) &= \left( \frac{-5\sin u \sin 2u + 4\cos u \cos 2u - 4\cos 2u}{\sqrt{2}}, -\frac{4}{\sqrt{2}}\sin 2u,- \frac{4}{\sqrt{2}}\cos 2u \right).
\end{align*}
Since $\gamma'(u) \times \gamma''(s)=(-1,\sin u (1+\cos^2 u\cos 2u), \sqrt{2}-(\cos u \sin^2 2u)/2)$, $\gamma$ is non-degenerate. 
Substituting $u=0$ and $u=\pi$ into these equations give 
\begin{align*}
\bt(0)= \bt(\pi) =  \left(\frac{1}{\sqrt{2}}, 0, \frac{1}{\sqrt{2}}\right), \ 
\bn(0) =\bn(\pi) =  \left(0, 1, 0\right),\ 
\bb(0) = \bb(\pi) = \left(-\frac{1}{\sqrt{3}}, 0, \frac{\sqrt{2}}{\sqrt{3}}\right).
\end{align*}  
Moreover, $t_{11}(0,\pi) = 1$, $t_{21}(0,\pi) = 0$ and $t_{33}(0, \pi)=1$, where $t_{ij}$ is the $(i,j)$ entry of the frame matrix between $(\gamma, \nu_1,\nu_2)$ and $(\gamma, \nu_1, \nu_2)$.
The curvature and torsion at $u=0$ and $u=\pi$ of the curve are 
\begin{eqnarray*}
(\kappa(0), \tau(0)) = (\sqrt{2}, -\sqrt{2}), \  (\kappa(\pi), \tau(\pi)) = (\sqrt{2},\sqrt{2}), 
\end{eqnarray*}
respectively. Thus, $\kappa(0) + \kappa(\pi)t_{11}(0, \pi) = 2\sqrt{2}$ and 
$\tau(0)\tau(\pi)(\kappa(0)^2 + \kappa(\pi)^2)t_{11}(0,\pi)-$ $\kappa(0)\kappa(\pi)$ $(\tau(0)^2+\tau(\pi)^2) = -16 < 0$. 
Moreover, the curvature $\alpha$ of the framed curve $(\gamma, \bn, \bb)$ is $$\alpha(u) = \sqrt{\frac{(\sin u \sin 2u + \cos 2u)^2 + 1}{2}}.$$ 
This satisfies $\alpha(0) = \alpha(\pi) = 1$ and $\alpha_u(0)=\alpha_u(\pi)=0$.
By Corollary \ref{self-translation-S1} (1) and Remark \ref{remark-S1assumption}, the translation surface $\bx^+$ has an $S_1^+$ singular point at $(0, \pi)$.
}
\end{example}
\begin{figure}[h!]
\begin{center}
\includegraphics[width=55mm]{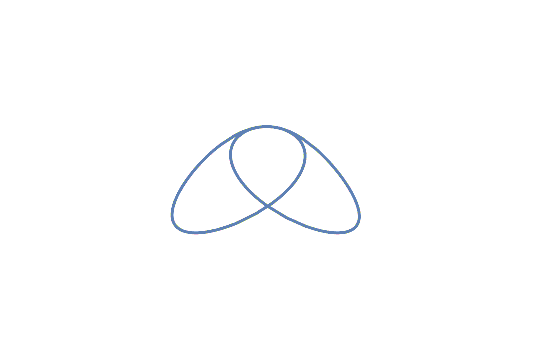}
\includegraphics[width=55mm]{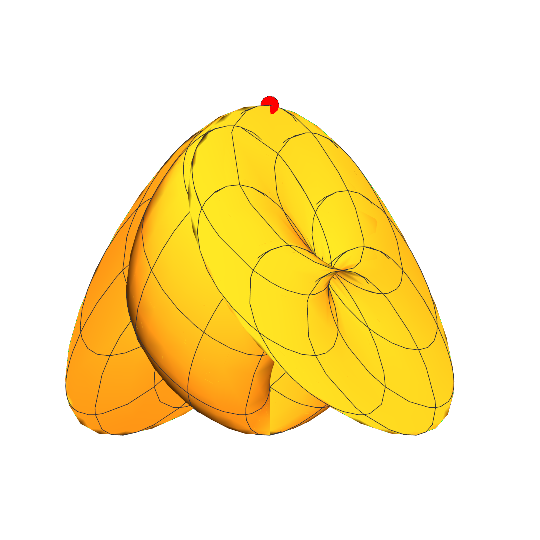}
\caption{The curve $\gamma$ and self translation surface $\bx^+$ in Example \ref{selfS1+exam}. $(0,\pi)$ is an $S_1^+$ singular point of $\bx^+$ and the dot in the figure on the right is the point $\bx^+ (0,\pi)$.}
\end{center}
\end{figure}

\begin{remark}{\rm 
For the curve in the above example, consider the self translation surface $\bx^-$. Then, $\kappa(0) - \kappa(\pi)t_{11}(0, \pi) = 0$ (note that the frame matrix between $(\gamma, \nu_1,\nu_2)$ and $(-\gamma, \nu_1, \nu_2)$ and the frame matrix between $(\gamma, \nu_1,\nu_2)$ and $(\gamma, \nu_1, \nu_2)$ are the same).
Thus, $\bx^-$ does not have an $S_1^+$ singular point at $(0,\pi)$.  This means that the analogue of Corollary \ref{self-translation-S1-} does not hold for an $S_1^+$ singular point.}
\end{remark}

\begin{remark}{\em
Using the framed curves $(\gamma,\nu_1,\nu_2)$ and $(\wgamma, \wnu_1,\wnu_2)$ in Example \ref{example-S1-}, 
we can construct a framed curve $(\gamma_c, \nu_{1c}, \nu_{2c})$ 
such that the self translation surface $\displaystyle \bx^{+}(u,v) = ({\gamma_c(u) + \gamma_c(v)})/{2}$ has an $S_1^-$ singular point at $(0,1)$. 
Intuitively, we cut off small neighbourhoods of $\gamma(0)$ and $\wgamma(0)$ from $\gamma$ and $\wgamma$, respectively, and then construct a curve by smoothly concatenate the resulting segments together, including their frames. More precisely, the construction can be carried out as follows.

Let $a, b$ be sufficiently small real numbers so that $a < 1-3b$. 
In the following, we define the framed curve $(\gamma_c, \nu_{1c}, \nu_{2c}) : [-a,1+b] \rightarrow \mathbb{R}^3 \times \Delta$ with the curvatures $(\ell_c,m_c,n_c,\alpha_c)$ by concatenating five framed curves which will be defined later
\begin{eqnarray}\label{curve}
(\gamma_c, \nu_{1c},\nu_{2c}) = 
\begin{cases}
(\gamma_{\varphi}, \nu_{1\varphi}, \nu_{2\varphi})  \ \text{on} \ [-a, a]\\
(\gamma_{[a,1-3b]}, \nu_{1[a,1-2b]}, \nu_{2[a,1-2b]}) \ \text{on} \ [a, 1-3b]\\
(\gamma_{[1-3b,1-2b]}, \nu_{1[1-3b,1-2b]},\nu_{2[1-3b,1-2b]}) \ \text{on} \ [1-3b, 1-2b]\\
(\gamma_{[1-2b,1-b]}, \nu_{1[1-2b,1-b]}, \nu_{2[1-2b,1-b]}) \ \text{on} \ [1-2b, 1-b]\\
(\wgamma_{\wvarphi}, \wnu_{1\wvarphi}, \wnu_{2\wvarphi}) \ \text{on} \ [1-b, 1+b]
\end{cases}
\end{eqnarray}

Let $\varphi : [-a,+a] \rightarrow \mathbb{R}$ be a $C^{\infty}$ function which satisfies
\begin{eqnarray*}
\varphi(u) = 
\begin{cases} 
u \ \ (u \in [- a/2,  a/2]) \\
0 \ \ (u \in [-a, a] \setminus ( - 2a/3,   2a/3)).
\end{cases}
\end{eqnarray*}
From the general theory of manifolds, such a function exists.
We define the map $(\gamma_{\varphi}, \nu_{1\varphi}, \nu_{2\varphi}) : [-a,a] \rightarrow \mathbb{R}^3 \times \Delta$ by 
\begin{eqnarray*}
(\gamma_{\varphi}, \nu_{1\varphi}, \nu_{2\varphi})(u) = (\gamma \circ \varphi, \nu_1 \circ \varphi, \nu_2 \circ \varphi)(u).
\end{eqnarray*}
By the properties of $\varphi$, we have $\gamma_\varphi(\pm a) = \gamma(0)$, $\gamma_{\varphi}(0) = \gamma(0)$ and $\gamma_{\varphi}(u) = \gamma(u)$ on $[- a/2,  a/2]$.
Moreover, 
\begin{eqnarray*}
\gamma_{\varphi, u}(u) = \varphi^{\prime}(u) \gamma_u(\varphi(u)) = \varphi^{\prime}(u) \alpha(\varphi(u))\mu(\varphi(u)).
\end{eqnarray*}
Thus, 
$(\gamma_{\varphi}, \nu_{1\varphi}, \nu_{2\varphi}) : [-a,a] \rightarrow \mathbb{R}^3 \times \Delta$ is a framed curve with the curvatures $(\ell_{\varphi}, m_{\varphi}, n_{\varphi}, \alpha_{\varphi})(u) = 
(\ell \circ \varphi, m \circ \varphi, n \circ \varphi, \varphi^{\prime}\alpha \circ \varphi)(u)$. 

Let $\wvarphi : [- b, b] \rightarrow \mathbb{R}$ be a $C^{\infty}$ function which satisfies the following properties: 
\begin{eqnarray*}
\wvarphi(u) = 
\begin{cases} 
u \ \ (u \in [- b/2,  b/2]) \\
0 \ \ (u \in [-b, b] \setminus (- 2b/3,  2b/3)).
\end{cases}
\end{eqnarray*}
We define the map $(\wgamma_{\wvarphi}, \wnu_{1\wvarphi}, \wnu_{2\wvarphi}) : [1-b,1+b] \rightarrow \mathbb{R}^3 \times \Delta$ by 
\begin{eqnarray*}
(\wgamma_{\wvarphi}, \wnu_{1\wvarphi}, \wnu_{2\wvarphi})(u) = (\wgamma \circ \wvarphi, \wnu_1 \circ \wvarphi, \wnu_2 \circ \wvarphi)(u - 1 ).
\end{eqnarray*}
By the properties of $\wvarphi$, we have 
$\wgamma_{\wvarphi}(1 \pm b) = \wgamma(0)$, 
$\wgamma_{\wvarphi}(1) = \wgamma(0)$, 
$\wgamma_{\wvarphi}(u) = \wgamma(u)$ on $[1 - b/2, 1 + b/2]$.
Moreover, 
\begin{eqnarray*}
\wgamma_{\wvarphi, u}(u) = \wvarphi^{\prime}(u -1 ) \wgamma_u(\varphi(u-1)) = \wvarphi^{\prime}(u-1) \walpha(\wvarphi(u-1))\wbmu(\wvarphi(u-1)).
\end{eqnarray*}
Thus, 
$(\wgamma_{\wvarphi}, \wnu_{1\wvarphi}, \wnu_{2\wvarphi}) : (1-b,1+b) \rightarrow \mathbb{R}^3 \times \Delta$ is a framed curve with the curvatures $(\well_{\wvarphi}, \wm_{\wvarphi}, \wn_{\wvarphi}, \walpha_{\wvarphi})(u)$ $=$ 
$(\well \circ \wvarphi, \wm \circ \wvarphi, \wn \circ \wvarphi, \wvarphi^{\prime}\walpha \circ \wvarphi)(u -1)$. 

Let $\varphi_{[1-3b,1-2b]} : \mathbb{R} \rightarrow [0, 1]$ be a $C^{\infty}$ function which satisfies 
\begin{eqnarray*}
\varphi_{[1-3b,1-2b]}(u) = 
\begin{cases} 
0 \ \ (u \in (-\infty, 1-3b]) \\
1 \ \ (u \in [1-2b, \infty)).
\end{cases}
\end{eqnarray*}

We define $\gamma_{[1-3b,1-2b]} : [1-3b, 1-2b] \rightarrow \mathbb{R}^2$ by 
$\gamma_{[1-3b,1-2b]}(u) = \gamma_\varphi(a) (1- \varphi_{[1-3b,1-2b]}(u)) +  \wgamma_{\wvarphi}(1-b)\varphi_{[1-3b,1-2b]}(u)$. Then $\gamma_{[1-3b,1-2b]} : [1-3b,1-2b] \rightarrow \mathbb{R}^3$ is a framed base curve (the image of $\gamma_{[1-3b,1-2b]}$ is a point or a part of line). We denote a framed curve which has $\gamma_{[1-3b,1-2b]}$ as the framed base curve $(\gamma_{[1-3b,1-2b]}, \nu_{1[1-3b,1-2b]},\nu_{2[1-3b,1-2b]})$.

Let $\gamma_{[a,1-3b]} : [a,1-3b] \rightarrow \mathbb{R}^3$ be the map $\gamma_{[a,1-3b]}(u) = \gamma_{\varphi}(a)$ for all $u \in [a,1-3b]$. We define a frame $(\nu_{1[a,1-3b]}, \nu_{2[a,1-3b]}, \mu_{[a,1-3b]})$ of $\gamma_{[a,1-3b]}$ such that it smoothly connects $(\nu_{1\varphi}(a), \nu_{2\varphi}(a), \mu_{\varphi}(a))$ and $(\nu_{1[1-3b,1-2b]}(1-3b), \nu_{2[1-3b,1-2b]}(1-3b), \mu_{[1-3b,1-2b]}(1-3b))$. See \cite{Fukunaga-Takahashi-2017} for more details.

Let $\gamma_{[1-2b,1-b]} : [1-2b,1-b] \rightarrow \mathbb{R}^3$ be the map $\gamma_{[1-2b,1-b]}(u) = \wgamma_{\wvarphi}(1-b)$ for all $u \in [1-2b,1-b]$. We define a frame $(\nu_{1[1-2b,1-b]}, \nu_{2[1-2b,1-b]}, \mu_{[1-2b,1-b]})$ of $\gamma_{[1-2b,1-b]}$ such that it smoothly connects $(\nu_{1[1-3b,1-2b]}(1-3b), \nu_{2[1-3b,1-2b]}(1-2b), \mu_{[1-3b,1-2b]}(1-3b))$ and 
$(\wnu_{1\wvarphi}(1-b), \wnu_{2\wvarphi}(1-b), \wmu_{\wvarphi}(1-b))$.

We define the framed curve $(\gamma_c, \nu_{1c}, \nu_{2c}) : [-a,1+b] \rightarrow \mathbb{R}^3 \times \Delta$ with the curvatures $(\ell_c,m_c,n_c,\alpha_c)$ by (\ref{curve}).
By the definitions of $\varphi$ and $\wvarphi$, $(\ell_c,m_c,n_c,\alpha_c)(u) = (\ell, m, n, \alpha)(u)$ around $u=0$ and $(\ell_c,m_c,n_c,\alpha_c)(u) = (\well,\wm, \wn, \walpha)(u-1)$ around $u=1$.
By Theorem \ref{S1} and Example \ref{example-S1-}, $(0 , 1) \in [-a, 1+b]$ is an $S_1^-$ singular point of the self translation surface $\displaystyle \bx^+(u,v) = ({\gamma_c(u)+ \gamma_c(v)})/{2}$ 

More generally, for framed curves $(\gamma, \nu_1, \nu_2)$ and $(\wgamma,\wnu_1,\wnu_2)$, the basic invariants of the generalised framed surface $(\bx, \nu_1, \nu_2)$ where $\displaystyle \bx(u,v) = ({\gamma(u)+ \wgamma(v)})/{2}$ determined by the curvature of the framed curves,
types of singularities of $\bx$ determined by the curvature of the framed curves. 
Thus, using the above method, if a singularity can be realized as a singular point of a translation surface, that the singularity can also be realized as a singular point of the self translation surface.
}\end{remark}

\section{Translation surfaces as framed surfaces}

We consider the case of the translation surface is a framed base surface under the tangency condition. 

Let $(\gamma,\nu_1,\nu_2) : I \rightarrow \mathbb{R}^3$ and $(\wgamma,\wnu_1,\wnu_2) : J \rightarrow \mathbb{R}^3$ be framed curves with curvatures $(\ell,m,n,\alpha)$ and $(\well,\wm,\wn,\walpha)$, respectively and $T$ be the frame matrix of $(\gamma,\nu_1,\nu_2)$ and $(\wgamma,\wnu_1,\wnu_2)$.
Then $(\bx,\nu_1,\nu_2):I \times J \to \R^3 \times \Delta$,
$$
\bx(u,v)=\gamma(u)+\wgamma(v), \ \nu_1(u,v)=\nu_1(u), \ \nu_2(u,v)=\nu_2(u)
$$
is a generalised framed surface, see \S 2. 
We give a condition that $\bx$ is a framed base surface (cf. Appendix \ref{GFS-FS}). 
In general theory, we need some assumptions by Proposition \ref{FBS.condition-translation}. 
However, in this case, we only need the condition that $t_{32}$ and $t_{31}$ are linearly dependent.

\begin{proposition}\label{translation-framed-surface}
Under the above notations, 
suppose that there exists $\theta : I \times J \rightarrow \mathbb{R}$ such that 
\begin{align}\label{condition-framed-surface}
t_{32}(u,v) \cos \theta(u,v) + t_{31}(u,v) \sin \theta(u,v)=0.
\end{align}
Then $(\bx,\bn, \mu) : I \times J \rightarrow \mathbb{R}^3 \times \Delta$ is a framed surface, where 
$$
\bn(u,v) = \sin \theta(u,v) \nu_1(u) + \cos \theta(u,v) \nu_2(u), \ \mu(u,v)=\nu_1(u) \times \nu_2(u)=\mu(u).
$$
\end{proposition}
\demo
By a direct calculation, we have 
\begin{align*}
&\bx_u(u,v)=\gamma_u(u)=\alpha(u)\mu(u), \ \bx_v(u,v)=\wgamma_v(v)=\walpha(v)\wmu(v),\\
&\bx_u(u,v)\times\bx_v(u,v)=\alpha(u)\walpha(v)\mu(u)\times(t_{31}(u,v)\nu_1(u)+t_{32}(u,v)\nu_2(u)+t_{33}(u,v)\mu(u))\\
&\textcolor{white}{\bx_u(u,v)\times\bx_v(u,v)}=\alpha(u)\walpha(v)(-t_{32}(u,v)\nu_1(u)+t_{31}(u,v)\nu_2(u)).
\end{align*}
If we take $\bn(u,v) = \sin \theta(u,v) \nu_1(u) + \cos \theta(u,v) \nu_2(u)$ and $\mu(u,v)=\nu_1(u) \times \nu_2(u)=\mu(u)$, then we have $\bx_u(u,v)\cdot\bn(u,v)=0$ and $\bx_v(u,v)\cdot\bn(u,v)=\alpha(u)\walpha(v)(t_{32}(u,v) \cos \theta(u,v) + t_{31}(u,v) \sin \theta(u,v))=0$. 
Therefore, $(\bx,\bn, \mu)$ is a framed surface.
\enD
By a direct calculation, we have the following.

\begin{proposition}\label{translation-framed-surface-basic-invariants}
Under the same assumptions as in Proposition \ref{translation-framed-surface}, 
the basic invariants of the framed surface $(\bx,\bn, \mu) : I \times J \rightarrow \mathbb{R}^3 \times \Delta$ are as follows.
\begin{align*}
&\begin{pmatrix}
a_1 & b_1 \\
a_2 & b_2
\end{pmatrix}(u,v)
=
\begin{pmatrix}
\alpha(u) & 0 \\
\walpha(v) t_{33}(u,v)  & 
\walpha(v)(-t_{32}(u,v)\sin\theta(u,v)+t_{31}(u,v)\cos\theta(u,v)) 
\end{pmatrix},\\
&\begin{pmatrix}
e_1 & f_1 & g_1 \\
e_2 & f_2 & g_2
\end{pmatrix}(u,v)\\
&=
\begin{pmatrix}
m(u)\sin\theta(u,v)+n(u)\cos\theta(u,v) & \theta_u(u,v)-\ell(u)& n(u)\sin\theta(u,v) - m(u)\cos\theta(u,v) \\
0  & \theta_v(u,v)& 0
\end{pmatrix}.
\end{align*}
\end{proposition}
\demo
We denote $\bt(u,v)=\bn(u,v)\times\mu(u,v)=\cos\theta(u,v)\nu_1(u)-\sin\theta(u,v)\nu_2(u)$. 
Since $\bx_u(u,v)= \ \alpha(u)\mu(u), \ \bx_v(u,v)=\walpha(v)\wmu(v),$ we have 
\begin{align*}
a_1(u,v)=& \ \bx_u(u,v)\cdot\mu(u)=\alpha(u)\mu(u)\cdot\mu(u)=\alpha(u),\\
a_2(u,v)=& \ \bx_v(u,v)\cdot\mu(u)=\walpha(v) \wmu(v) \cdot\mu(u)
= \walpha(v)t_{33}(u,v),\\
b_1(u,v)=& \ \bx_u(u,v)\cdot\bt(u,v)=\alpha(u)\mu(u)\cdot(\cos\theta(u,v)\nu_1(u)-\sin\theta(u,v)\nu_2(u))=0,\\
b_2(u,v)=& \ \bx_v(u,v)\cdot\bt(u,v)=\walpha(v)(-t_{32}(u,v)\sin\theta(u,v)+t_{31}(u,v)\cos\theta(u,v)).\end{align*}
Moreover, 
\begin{align*}
\bn_u(u,v)=& \ \theta_u(u,v)\cos\theta(u,v)\nu_1(u)+\sin\theta(u,v)(\ell(u)\nu_2(u)+m(u)\mu(u))\\
&-\theta_u(u,v)\sin\theta(u,v)\nu_2(u)+\cos\theta(u,v)(-\ell(u)\nu_1(u)+n(u)\mu(u))\\
=& \ (\theta_u(u,v)-\ell(u))(\cos\theta(u,v)\nu_1(u)-\sin\theta(u,v)\nu_2(u))\\
&+(m(u)\sin\theta(u,v)+n(u)\cos\theta(u,v)),\\
\bn_v(u,v)=& \ \theta_v(u,v)(\cos\theta(u,v)\nu_1(u)-\sin\theta(u,v)\nu_2(u)) ,\\
\bt_u(u,v)=& \ -(\theta_u(u,v)-\ell(u))(\cos\theta(u,v)\nu_2(u)-\sin\theta(u,v)\nu_1(u))\\
&+(m(u)\cos\theta(u,v)-n(u)\sin\theta(u,v))\mu(u),\\
\bt_v(u,v)=& \ -\theta_v(u,v)(\cos\theta(u,v)\nu_2(u)+\sin\theta(u,v)\nu_1(u)),
\end{align*}
we have 
\begin{align*}
e_1(u,v)=& \ \bn_u(u,v)\cdot\mu(u)=m(u)\sin\theta(u,v)+n(u)\cos\theta(u,v),\\
f_1(u,v)=& \ \bn_u(u,v)\cdot\bt(u,v)=\theta_u(u,v)-\ell(u),\\
e_2(u,v)=& \ \bn_v(u,v)\cdot\mu(u)=0,\\
f_2(u,v)=& \ \bn_v(u,v)\cdot\bt(u,v)=\theta_v(u,v),\\
g_1(u,v)=& \ -\bt_u(u,v)\cdot\mu(u)=n(u)\sin\theta(u,v)-m(u)\cos\theta(u,v),\\
g_2(u,v)=& \ -\bt_v(u,v)\cdot\mu(u)=0.
\end{align*}
\enD
By using Proposition \ref{translation-framed-surface-basic-invariants}, we have the curvature $C^F=(J^F,K^F,H^F)$ of the framed surface $(\bx,\bn, \mu)$, see Appendix \ref{GFS-FS}.
\begin{align*}
J^F(u,v) &= \alpha(u)\walpha(v) \left(-t_{32}(u,v)\sin\theta(u,v)+t_{31}(u,v)\cos\theta(u,v) \right),\\
K^F(u,v) &= \theta_v(u,v) \left(m(u)\sin\theta(u,v)+n(u)\cos\theta(u,v)\right), \\
H^F(u,v) &= -\frac{1}{2}\left\{ \alpha(u)\theta_v(u,v)-\walpha(v)t_{33}(u,v)(\theta_u(u,v)-\ell(u))-\walpha(v)(m(u)\sin\theta(u,v) + \right.\\
&\quad \left. n(u)\cos\theta(u,v))(-t_{32}(u,v)\sin\theta(u,v)+t_{31}(u,v)\cos\theta(u,v)) \right\}.
\end{align*}

We investigate criteria for singular points of translation surfaces as framed base surfaces.
We use the following lemma.
\begin{lemma}\label{lemma-framed-surface}
Let $(\gamma, \nu_1, \nu_2) : I \rightarrow \mathbb{R}^3 \times \Delta$ and $(\wgamma, \wnu_1, \wnu_2) : J \rightarrow \mathbb{R}^3 \times \Delta$ be framed curves with curvatures $(\ell, m, n, \alpha)$ and $(\well, \wm, \wn, \walpha)$, respectively. 
Suppose that there exists $\theta : I \times J \rightarrow \mathbb{R}$ such that $t_{32}(u,v) \cos \theta(u,v) + t_{31}(u,v) \sin \theta(u,v)=0$.
Let $p_0 = (u_0, v_0)$ be a point with $\mu(u_0) \times \wmu(v_0) = 0$. The following relational equations hold.
\par
$(1)$ $m(u_0)\sin\theta(p_0) + n(u_0)\cos\theta(p_0)=0$.
\par
$(2)$ $-(\wm(v_0)t_{12}(p_0)+\wn(v_0)t_{22}(p_0))\cos\theta(p_0) -(\wm(v_0)t_{11}(p_0)+\wn(v_0)t_{21}(p_0))\sin\theta(p_0)=0$.
\par
$(3)$ \begin{align*}
&(\ell(u_0)n(u_0)+m_u(u_0)-2n(u_0)\theta_u(p_0))\sin\theta(p_0) \\
&\quad +(-\ell(u_0)m(u_0)+n_u(u_0)+2m(u_0)\theta_u(p_0))\cos\theta(p_0)=0.
\end{align*}
\par
$(4)$ \begin{align*}
&\left\{ (\ell(u_0) - \theta_u(p_0))(\wm(v_0)t_{11}(p_0)+\wn(v_0)t_{21}(p_0))+\theta_v(p_0)m(u_0)t_{33}(p_0) \right\}\cos\theta(p_0)\\ 
& \quad -\left\{ (\ell(u_0) - \theta_u(p_0))(\wm(v_0)t_{12}(p_0)+\wn(v_0)t_{22}(p_0))+\theta_v(p_0)n(u_0)t_{33}(p_0) \right\}\sin\theta(p_0)=0.
\end{align*}
\par
$(5)$ \begin{align*}
&\bigl\{ (-\wm_v(v_0) +\well(v_0)\wn(v_0))t_{12}(p_0)-(\well(v_0)\wm(v_0) +\wn_v(v_0))t_{22}(p_0) \\
& -2(\wm(v_0)t_{11}(p_0) +\wn(v_0)t_{21}(p_0))\theta_v(p_0)\bigr\}\cos\theta(p_0) \\ 
& +\bigl\{ (-\wm_v(v_0) +\well(v_0)\wn(v_0))t_{11}(p_0)-(\well(v_0)\wm(v_0) +\wn_v(v_0))t_{21}(p_0) \\ 
& +2(\wm(v_0)t_{12}(p_0) +\wn(v_0)t_{22}(p_0))\theta_v(p_0)\bigr\}\sin\theta(p_0)=0.
\end{align*}
\par
$(6)$ \begin{align*}
&\bigl\{3(\theta_{uu}(p_0)m(u_0)+\theta_{u}(p_0)m_u(u_0))-\ell(u_0)m_u(u_0)+n_{uu}(u_0)-2\ell_u(u_0)m(u_0) \bigr\}\cos\theta(p_0) \\ 
& - \bigl\{3(\theta_{uu}(p_0)n(u_0)+\theta_{u}(p_0)n_u(u_0))-\ell(u_0)n_u(u_0)-m_{uu}(u_0)-2\ell_u(u_0)n(u_0) \bigr\}\sin\theta(p_0)=0.
\end{align*}
\par
$(7)$ \begin{align*}
&\bigl\{(\ell_u(u_0)+m(u_0)n(u_0)-\theta_{uu}(p_0))(\wm(v_0)t_{11}(p_0)+\wn(v_0)t_{21}(p_0))\\
&+(\theta_v(p_0)(\ell(u_0)n(u_0)+m_u(u_0))+2\theta_{uv}(p_0)m(u_0))t_{33}(p_0)\bigr\}\cos\theta(p_0) \\ 
& \quad -\bigl\{(\ell_u(u_0)-m(u_0)n(u_0)-\theta_{uu}(p_0))(\wm(v_0)t_{12}(p_0)+\wn(v_0)t_{22}(p_0))\\
& \quad +(\theta_v(p_0)(-\ell(u_0)m(u_0)+n_u(u_0))+2\theta_{uv}(p_0)n(u_0))t_{33}(p_0)\bigr\}\sin\theta(p_0)=0.
\end{align*}
\par
$(8)$ \begin{align*}
&\bigl\{(2\theta_{uv}(p_0)\wm(v_0)-(\well(v_0)\wn(v_0)-\wm_v(v_0))(\theta_u(p_0)-\ell(u_0)))t_{12}(p_0)\\
&+(2\theta_{uv}(p_0)\wn(v_0)+(\well(v_0)\wm(v_0)+\wn_v(v_0))(\theta_u(p_0)-\ell(u_0)))t_{22}(p_0)-\theta_{vv}(p_0)n(u_0)t_{33}(p_0)\bigr\}\cos\theta(p_0) \\ 
& \quad -\bigl\{(2\theta_{uv}(p_0)\wm(v_0)-(\well(v_0)\wn(v_0)-\wm_v(v_0))(\theta_u(p_0)-\ell(u_0)))t_{11}(p_0)\\
& \quad +(2\theta_{uv}(p_0)\wn(v_0)+(\well(v_0)\wm(v_0)+\wn_v(v_0))(\theta_u(p_0)-\ell(u_0)))t_{21}(p_0)\\
& \quad -\theta_{vv}(p_0)m(u_0)t_{33}(p_0)\bigr\}\sin\theta(p_0)=0.
\end{align*}
\par
$(9)$ \begin{align*}
&\bigl\{3(\theta_v(p_0)(\well(v_0)\wn(v_0)-\wm_v(v_0))-\theta_{vv}(p_0)\wm(v_0))t_{11}(p_0)\\
&+(\well_v(v_0)\wn(v_0)-\wm_{vv}(v_0)+2\well(v_0)\wn_v(v_0))t_{12}(p_0)\\
&-3(\theta_v(p_0)(\well(v_0)\wm(v_0)+\wn_v(v_0))+\theta_{vv}(p_0)\wn(v_0))t_{21}(p_0)\\
&-(\well_v(v_0)\wm(v_0)+\wn_{vv}(v_0)+2\well(v_0)\wm_v(v_0))t_{22}(p_0)\bigr\}\cos\theta(p_0) \\ 
& \quad +\bigl\{(\well_v(v_0)\wn(v_0)-\wm_{vv}(v_0)+2\well(v_0)\wn_v(v_0))t_{11}(p_0)\\
& \quad -3(\theta_v(p_0)(\well(v_0)\wn(v_0)-\wm_v(v_0))-\theta_{vv}(p_0)\wm(v_0))t_{12}(p_0)\\
& \quad -(\well_v(v_0)\wm(v_0)+\wn_{vv}(v_0)+2\well(v_0)\wm_v(v_0))t_{21}(p_0)\\
& \quad +3(\theta_v(p_0)(\well(v_0)\wm(v_0)+\wn_v(v_0))+\theta_{vv}(p_0)\wn(v_0))t_{22}(p_0)\bigr\}\sin\theta(p_0)=0.
\end{align*}
\end{lemma}
\demo
By differentiating $t_{32}(u,v) \cos \theta(u,v) + t_{31}(u,v) \sin \theta(u,v)=0$ from the first to the third order, we obtain the above equations.
\enD
\begin{corollary}\label{rel_theta}
Under the same assumptions as in Lemma \ref{lemma-framed-surface}, 
suppose that $\gamma : I \rightarrow \mathbb{R}^3$ and $\wgamma : J \rightarrow \mathbb{R}^3$ are non-degenerate curves with arc-length parameter. 
\par
$(1)$ $\sin\theta(p_0) = 0, \cos\theta(p_0)=\pm 1$.
\par
$(2)$ $t_{12}(p_0) = 0$.
\par
$(3)$ $\displaystyle \theta_u(p_0) = {\tau(u_0)}/{2}$.
\par
$(4)$ $\displaystyle \wkappa(v_0)(\tau(u_0) - \theta_u(p_0))t_{11}(p_0) + \theta_v(p_0)\kappa(u_0)t_{33}(p_0)= 0$.
\par
$(5)$ $\displaystyle \wtau(v_0)t_{22}(p_0)+2t_{11}(p_0) \theta_v(p_0)=0$.
\par
$(6)$ $\displaystyle 3\theta_{uu}(p_0)\kappa(u_0)+\frac{1}{2}\tau(u_0)\kappa_u(u_0)-2\tau_u(u_0)\kappa(u_0)=0$.
\par
$(7)$ $\displaystyle \wkappa(v_0)t_{11}(p_0)(\theta_{uu}(p_0)-\tau_u(u_0))-(\theta_v(p_0)\kappa_u(u_0)+2\theta_{uv}(p_0)\kappa(u_0))t_{33}(p_0)=0$.
\par
$(8)$ $\displaystyle (2\theta_{uv}(p_0)\wkappa(v_0)+\wkappa_v(v_0)(\theta_{u}(p_0)-\tau(u_0)))t_{11}(p_0)-\theta_{vv}(p_0)\kappa(u_0)t_{33}(p_0)=0$.
\par
$(9)$ $\displaystyle 3(\theta_v(p_0)\wkappa_v(v_0)+\theta_{vv}(p_0)\wkappa(v_0))t_{11}(p_0)+(\wtau_v(v_0)\wkappa(v_0)+2\wtau(v_0)\wkappa_v(v_0))t_{22}(p_0)=0$.
\end{corollary}
\demo
Since $\gamma$ and $\wgamma$ are non-degenerate curves with arc-length parameter, then the curvatures of $\gamma$ and $\wgamma$ are given by $(\ell,m,n,\alpha)=(\tau,-\kappa,0,1)$ and $(\well,\wm,\wn,\walpha)=(\wtau,-\wkappa,0,1)$. By calculating $(\ell,m,n,\alpha)$ and $(\well,\wm,\wn,\walpha)$ in the equation in Lemma \ref{lemma-framed-surface} as $(\tau,-\kappa,0,1)$ and $(\wtau,-\wkappa,0,1)$, we obtain the above equations.
\enD

The discriminant function (or, the signed area density function) 
$\lambda : I \times J \rightarrow \mathbb{R}$ is given by 
\begin{align*}
\lambda(u,v) &= \det(\bx_u, \bx_v,\bn)(u,v) \\
&= \left((\bx_u \times \bx_v) \cdot \bn\right)(u,v) \\
&= \alpha(u)\walpha(v)\left( -t_{32}(u,v) \sin\theta(u,v) + t_{31}(u,v)\cos\theta(u,v) \right).
\end{align*}

We divide into four cases of singular points of the translation framed surface under the linearly dependent condition:
\begin{itemize}
\item[(I)] $\alpha(u_0) \neq 0, \ \walpha(v_0) \neq 0, \ \mu(u_0) \times \wmu(v_0) = 0$,
\item[(II)] $\alpha(u_0) = 0, \ \walpha(v_0) \neq 0, \ \mu(u_0) \times \wmu(v_0) = 0$,
\item[(III)] $\alpha(u_0) \neq 0, \ \walpha(v_0) = 0, \ \mu(u_0) \times \wmu(v_0) = 0$,
\item[(IV)] $\alpha(u_0) = 0, \ \walpha(v_0) = 0, \ \mu(u_0) \times \wmu(v_0) = 0$.
\end{itemize}
\subsection{$\rm(I)$ The case of $\alpha(u_0) \neq 0, \walpha(v_0) \neq 0, \mu(u_0) \times \wmu(v_0) = 0$}\label{case1}

We may take the discriminant function (or, the signed area density function) as   $\Lambda : I \times J \rightarrow \mathbb{R}$, 
$$
\Lambda(u,v) =  -t_{32}(u,v) \sin\theta(u,v) + t_{31}(u,v)\cos\theta(u,v)
$$ 
(cf. Appendix \ref{criterion}). 
By a direct calculation and Lemma \ref{lemma-framed-surface}, we have
\begin{align*}
\Lambda_u(u,v) &= (\ell(u)t_{31}(u,v)-n(u)t_{33}(u,v))\sin\theta(u,v)\\ 
&\quad + (\ell(u)t_{32}(u,v)+m(u)t_{33}(u,v))\cos\theta(u,v), \\
\Lambda_v(u,v) &= (\wm(v)t_{12}(u,v)+\wn(v)t_{22}(u,v))\sin\theta(u,v)\\ 
&\quad - (\wm(v)t_{11}(u,v)+\wn(v)t_{21}(u,v))\cos\theta(u,v) 
\end{align*}
and 
\begin{align*}
\Lambda_u(p_0) &= t_{33}(p_0)(-n(u_0)\sin\theta(p_0)+m(u_0)\cos\theta(p_0)),\\
\Lambda_v(p_0) &= (\wm(v_0)t_{12}(p_0)+\wn(v_0)t_{22}(p_0))\sin\theta(p_0)\\
&\quad - (\wm(v_0)t_{11}(p_0)+\wn(v_0)t_{21}(p_0))\cos\theta(p_0).
\end{align*}
Moreover, 
\begin{align*}
\Lambda_{uu}(p_0) &= t_{33}(p_0)(-n_{u}(u_0)\sin\theta(p_0)+m_{u}(u_0)\cos\theta(p_0)),\\
\Lambda_{uv}(p_0) &= 0, \\
\Lambda_{vv}(p_0) &= ((\wm_{v}(v_0)-\well(v_0)\wn(v_0))t_{12}(p_0)+(\wn_{v}(v_0)+\well(v_0)\wm(v_0))t_{22}(p_0))\sin\theta(p_0)\\
&\quad - ((\wm_v(v_0)-\well(v_0)\wn(v_0))t_{11}(p_0)+(\wn_v(v_0)+\well(v_0)\wm(v_0))t_{21}(p_0))\cos\theta(p_0).
\end{align*}

\begin{theorem}\label{criteria-case-I}
Under the assumptions $\rm (I)$, we have the following.
\par
$(A)$ Suppose that $p_0$ is a non-degenerate singular point of $\bx$. 
\par
$(1)$ $p_0$ is a cuspidal edge of $\bx$ if and only if 
\begin{align*}
&\alpha(u_0)\theta_v(p_0)-\walpha(v_0)t_{33}(p_0)(\theta_u(p_0)-\ell(u_0))\neq 0,\\
&\walpha(v_0)(-n(u_0)\sin\theta(p_0)+m(u_0)\cos\theta(p_0))-\\
&\alpha(u_0)\left\{(\wm(v_0)t_{12}(p_0)+\wn(v_0)t_{22}(p_0))\sin\theta(p_0)-(\wm(v_0)t_{11}(p_0)+\wn(v_0)t_{21}(p_0))\cos\theta(p_0)\right\} \neq 0.
\end{align*}
\par  
$(2)$ Suppose that $\alpha_u(u)=0$ and $\walpha_v(v)=0$ for all $(u,v) \in I \times J$. 
Then $p_0$ is a swallowtail of $\bx$ if and only if 
\begin{align*}
&\alpha(u_0)\theta_v(p_0)-\walpha(v_0)t_{33}(p_0)(\theta_u(p_0)-\ell(u_0))\neq 0, \\
&\walpha(v_0)(-n(u_0)\sin\theta(p_0)+m(u_0)\cos\theta(p_0))-\\
&\alpha(u_0)\left\{(\wm(v_0)t_{12}(p_0)+\wn(v_0)t_{22}(p_0))\sin\theta(p_0)-(\wm(v_0)t_{11}(p_0)+\wn(v_0)t_{21}(p_0))\cos\theta(p_0)\right\} = 0,\\
&((\wm_v(v_0)-\well(v_0)\wn(u_0)) t_{12}(p_0)+(\wn_v(v_0)+\well(v_0)\wm(u_0)) t_{22}(p_0)-n_u(u_0)t_{33}(p_0))\sin \theta(p_0)\\
&-((\wm_v(v_0)-\well(v_0)\wn(u_0)) t_{11}(p_0)+(\wn_v(v_0)+\well(v_0)\wm(u_0)) t_{21}(p_0)+m_u(u_0)t_{33}(p_0))\cos \theta(p_0) \not=0.
\end{align*} 
\par
$(3)$ Suppose that $\alpha_u(u)=0$ and $\walpha_v(v)=0$ for all $(u,v) \in I \times J$. 
Then $p_0$ is a cuspidal cross cap of $\bx$ if and only if 
\begin{align*}
&\walpha(v_0)(-n(u_0)\sin\theta(p_0)+m(u_0)\cos\theta(p_0))-\\
&\alpha(u_0)((\wm(v_0)t_{12}(p_0)+\wn(v_0)t_{22}(p_0))\sin\theta(p_0)-(\wm(v_0)t_{11}(p_0)+\wn(v_0)t_{21}(p_0))\cos\theta(p_0)) \neq 0,\\
&\alpha(u_0)\theta_v(p_0)-\walpha(v_0)t_{33}(p_0)(\theta_u(p_0)-\ell(u_0))= 0, 
\end{align*}
\begin{align*}
&\alpha(u_0)(\theta_{uv}(u_0,v_0)((\wm(v_0)t_{12}(p_0)+\wn(v_0)t_{22}(p_0))\sin\theta(p_0)\\
&-(\wm(v_0)t_{11}(p_0)+\wn(v_0)t_{21}(p_0))\cos\theta(p_0)) \\
&-\theta_{vv}(p_0) t_{33}(p_0)(-n(u_0)\sin \theta(p_0)+m(u_0)\cos \theta(p_0))) \\
&-\walpha(v_0) t_{33}(p_0)(\theta_{uu}(p_0) ((\wm(v_0)t_{12}(p_0)+\wn(v_0)t_{22}(p_0))\sin\theta(p_0)\\
&-(\wm(v_0)t_{11}(p_0)+\wn(v_0)t_{21}(p_0))\cos\theta(p_0)) \\
&-\theta_{uv}(p_0) t_{33}(p_0)(-n(u_0)\sin\theta(p_0)+m(u_0)\cos\theta(p_0))\\
&-\ell_u(u_0)((\wm(v_0)t_{12}(p_0)+\wn(v_0)t_{22}(p_0))\sin\theta(p_0)\\
&-(\wm(v_0)t_{11}(p_0)+\wn(v_0)t_{21}(p_0))\cos\theta(p_0))) \neq 0.
\end{align*}
\par 
\par
$(B)$ Suppose that $p_0$ is a degenerate singular point of $\bx$, that is, $\Lambda_u(p_0)=\Lambda_v(p_0)=0$. 
Then $p_0$ is neither cuspidal lips nor cuspidal beaks.
\end{theorem}
\demo
$(A)$ Suppose that $\Lambda_u(p_0) \not=0$. 
In the case of $\Lambda_v(p_0) \not=0$, we can also prove by the similar calculations, we omit it.
\par
$(1)$ There exists $u:J \to \R$ such that $\Lambda(u(v),v)=0$ around $p_0$. 
Since the singular set of $\bx$ is $S(\bx)=\{(u(v),v) \in I \times J\}$, the singular curve is given by $\delta:(J,v_0) \to (I \times J,p_0)$. 
Since $J^F(u(t),t)=K^F(u(t),t)=0$ and 
\begin{align*}
H^F(u(t),t) &= -\frac{1}{2} \left\{ \alpha(u(t))\theta_v(u(t),t)-\widetilde{\alpha}(t)t_{33}(u(t),t)(\theta_u(u(t),t)-\ell(u(t)))\right\},
\end{align*}
$\bx$ is a front if and only if $\alpha(u_0)\theta_v(p_0)-\widetilde{\alpha}(v_0)t_{33}(p_0)(\theta_u(p_0)-\ell(u_0)) \not=0$ by Proposition \ref{Legendre.immersion}. 
By $d\bx |_{\delta(t)}=\alpha(u(t)) \mu(u(t)) du+\widetilde{\alpha}(t)t_{33}(u(t),t)\mu(u(t)) dv$, 
we can take a null vector field $\eta(u,v)=-\widetilde{\alpha}(v)t_{33}(u,v) \partial/\partial u +\alpha(u) \partial/\partial v$. 
It follows that 
\begin{align*}
\eta \Lambda (p_0) =& -\walpha(v_0)(-n(u_0)\sin\theta(p_0)+m(u_0)\cos\theta(p_0))+\\
&\alpha(u_0)\left\{(\wm(v_0)t_{12}(p_0)+\wn(v_0)t_{22}(p_0))\sin\theta(p_0)+(\wm(v_0)t_{11}(p_0)+\wn(v_0)t_{21}(p_0))\cos\theta(p_0)\right\}.
\end{align*} 
By the criterion of the cuspidal edge, we have the assertion. 
\par
$(2)$ By using Lemma \ref{lemma-framed-surface} and a direct calculation, we have 
$\eta \eta \Lambda (p_0)= \Lambda_{uu}(p_0)+\Lambda_{vv}(p_0)$ at $p_0$. 
By the criterion of the swallowtail, we have the assertion. 
\par
$(3)$ By a direct calculation, we have 
\begin{align*}
(\bx \circ \delta)^{\prime}(t) &= ( \alpha(t) u^{\prime}(t)+ \walpha(t) t_{33}(u(t),t) )\mu(t), \\
\bn(\delta(t)) &= \sin\theta(\delta(t))\nu_1(\delta(t)) + \cos\theta(\delta(t))\nu_2(\delta(t)), \\
d\bn(\eta) &= \left\{ \walpha(t)t_{33}(\delta(t))(\theta_u(\delta(t))-\ell(u(t)))-\alpha(u(t))\theta_v(\delta(t)) \right\}\cos\theta(\delta(t))\nu_1(u(t)) \\
&\quad + \left\{ \walpha(t)t_{33}(\delta(t))(-\theta_u(\delta(t))+\ell(u(t)))+\alpha(u(t))\theta_v(\delta(t)) \right\}\sin\theta(\delta(t))\nu_2(u(t))\\
&\quad +\walpha(t)t_{33}(\delta(t))(m(u(t))\sin\theta(\delta(t))+m(u(t))\cos\theta(\delta(t)))\mu(u(t)).
\end{align*}
Therefore 
\begin{align*}
\phi_{\bx}(t) &= \det((\bx \circ \delta)^{\prime}, \bn \circ \delta, d\bn(\eta))(t) \\
&= (\alpha(u(t)) u^{\prime}(t) + \walpha(t)t_{33}(\delta(t)))
\left\{ \walpha(t)t_{33}(\delta(t))(\theta_u(\delta(t))-\ell(u(t))) - 
\alpha(u(t))\theta_v(\delta(t)) \right\}.
\end{align*}
By $\phi_{\bx}(t_0)=0, \phi'_{\bx}(t_0) \not=0$ and the criterion of the cuspidal cross cap, we have the assertion. 
\par
$(B)$ By the assumption, Lemma \ref{lemma-framed-surface} (1) and (2), we have $m(u_0)=n(u_0)=\wm(v_0)=\wn(v_0)=0$. 
By a direct calculation, we have 
\begin{align*}
&\det(\text{\rm Hess} \lambda(p)) 
= \det\begin{pmatrix} 
\Lambda_{uu}(p_0) & \Lambda_{uv}(p_0) \\
\Lambda_{vu}(p_0) & \Lambda_{vv}(p_0)
\end{pmatrix} \\
&= t_{33}(p_0) (-n_u(u_0)\sin \theta(p_0)+m_u(u_0)\cos \theta(p_0))\\
& \quad \left\{(\wm_v(v_0)t_{12}(p_0)+\wn_v(v_0)t_{22}(p_0))\sin\theta(p_0)-(\wm_v(v_0)t_{11}(p_0)+\wn_v(v_0)t_{21}(p_0))\cos\theta(p_0) \right\}.
\end{align*}
If $\det(\text{\rm Hess} \lambda(p)) \not=0$, then $\theta_u(p_0)-\ell(u_0)=0$ and $\theta_v(p_0)=0$ using Lemma \ref{lemma-framed-surface} $(7)$ and $(8)$. 
It follows that $\bx$ is not a front. 
Hence, $p_0$ is neither cuspidal lips  nor cuspidal beaks.
\enD

\begin{corollary}\label{criteria-case-I-Frenet}
Under the assumptions $\rm (I)$, suppose that $\gamma : I \rightarrow \mathbb{R}^3$ and $\wgamma : J \rightarrow \mathbb{R}^3$ are non-degenerate curves with arc-length parameter and $p_0$ is a non-degenerate singular point of $\bx$.
\par
$(1)$ $p_0$ is a cuspidal edge of $\bx$ if and only if 
\begin{align*}
&\tau(u_0)(\kappa(u_0)-\wkappa(v_0)t_{11}(p_0)) \not=0, \  \kappa(u_0)+\wkappa(v_0)t_{11}(p_0) \not=0.
\end{align*}
\par  
$(2)$  $p_0$ is a swallowtail of $\bx$ if and only if 
\begin{align*}
\tau(u_0)(\kappa(u_0)-\wkappa(v_0)t_{11}(p_0)) \not=0,  \kappa(u_0)+\wkappa(v_0)t_{11}(p_0) =0, 
-\kappa_u(u_0) t_{33}(p_0)+\wkappa_v(v_0) t_{11}(p_0) \not=0.
\end{align*}
\par
$(3)$ $p_0$ is a cuspidal cross cap of $\bx$ if and only if 
\begin{align*}
\tau(u_0)(\kappa(u_0)-\wkappa(v_0)t_{11}(p_0)) =0, 
\kappa(u_0)+\wkappa(v_0)t_{11}(p_0) \not= 0, 
-\kappa_u(u_0) t_{33}(p_0)+\wkappa_v(v_0) t_{11}(p_0) \not=0.
\end{align*}
\end{corollary}
\demo
Since $\gamma$ and $\wgamma$ are non-degenerate curves with arc-length parameter, then the curvatures of $\gamma$ and $\wgamma$ are given by $(\ell,m,n,\alpha)=(\tau,-\kappa,0,1)$ and $(\well,\wm,\wn,\walpha)=(\wtau,-\wkappa,0,1)$. 
By calculating $(\ell,m,n,\alpha)$ and $(\well,\wm,\wn,\walpha)$ in the equations in Theorem  \ref{criteria-case-I} as $(\tau,-\kappa,0,1)$ and $(\wtau,-\wkappa,0,1)$ and using Corollary \ref{rel_theta}, we obtain the results.
\enD

\begin{remark}{\rm
If $\alpha(u_0) =\walpha(v_0) = 0$ and $\mu(u_0) \times \wmu(v_0) \not= 0$, then 
the translation framed base surface $\bx$ never have swallowtail singular points \cite{Fukunaga-Takahashi-2022}.
}
\end{remark}


\subsection{$\rm(II)$ The case of $\alpha(u_0) = 0, \walpha(v_0) \neq 0, \mu(u_0) \times \wmu(v_0) = 0$}\label{case2}

By a direct calculation, we have 
\begin{align*}
\lambda(u,v)&= \alpha(u)\walpha(v)\left( -t_{32}(u,v) \sin\theta(u,v) + t_{31}(u,v)\cos\theta(u,v) \right), \\
\lambda_u(u,v)&= \alpha_u(u)\walpha(v)\left( -t_{32}(u,v) \sin\theta(u,v) + t_{31}(u,v)\cos\theta(u,v) \right)\\
& \quad +\alpha(u)\walpha(v)\left( -t_{32}(u,v) \sin\theta(u,v) + t_{31}(u,v)\cos\theta(u,v) \right)_u,\\ 
\lambda_v(u,v)&= \alpha(u)\walpha_v(v)\left( -t_{32}(u,v) \sin\theta(u,v) + t_{31}(u,v)\cos\theta(u,v) \right)\\
& \quad +\alpha(u)\walpha(v)\left( -t_{32}(u,v) \sin\theta(u,v) + t_{31}(u,v)\cos\theta(u,v) \right)_v.
\end{align*}
By $\alpha(u_0) = 0$, $t_{32}(p_0) = t_{31}(p_0) = 0$, we have $\lambda_u(p_0) = \lambda_v(p_0) = 0$. 
It follows that $p_0$ is a degenerate singular point of $\bx$ and rank ${d\bx}=1$ at $p_0$. 

\begin{theorem}\label{criteria-case-II}
Under the assumptions $\rm (II)$, we have the following.
\par
$(1)$ $p_0$ is never a cuspidal lips of $\bx$. 
\par
$(2)$ $p_0$ is a cuspidal beaks of $\bx$ if and only if 
$$
\theta_u(p_0) - \ell(u_0) \neq 0,\ \alpha_u(u_0) \neq 0,\ -n(u_0)\sin \theta(p_0) + m(u_0) \cos \theta(p_0) \neq 0
$$
and 
$\left(\wm(v_0)t_{12}(p_0)+ \wn(v_0)t_{22}(p_0)\right)\sin\theta(p_0) -
\left(\wm(v_0)t_{11}(p_0) + \wn(v_0)t_{21}(p_0)\right)\cos\theta(p_0) \not=0$.
\end{theorem}
\demo
$(1)$ By a direct calculation, we have 
\begin{align*}
\lambda_{uu}(p_0)&=2\alpha_u(u_0)\walpha(v_0)t_{33}(p_0)(-n(u_0)\sin\theta(p_0)+m(u_0)\cos\theta(p_0)),\\
\lambda_{uv}(p_0)&=\alpha_u(u_0)\walpha(v_0)\bigl((\wm(v_0)t_{12}(p_0)+ \wn(v_0)t_{22}(p_0))\sin\theta(p_0)\\
&\quad-(\wm(v_0)t_{11}(p_0) + \wn(v_0)t_{21}(p_0))\cos\theta(p_0)\bigr),\\
\lambda_{vv}(p_0)&=0.
\end{align*}
Therefore, we have 
\begin{align*}
\det (\text{\rm Hess} \lambda(p_0)) &=-\alpha_u(u_0)^2\walpha(v_0)^2\bigl((\wm(v_0)t_{12}(p_0)+ \wn(v_0)t_{22}(p_0)) \sin\theta(p_0) \\
& \quad -(\wm(v_0)t_{11}(p_0) + \wn(v_0)t_{21}(p_0))\cos\theta(p_0)\bigr)^2
\end{align*}
and hence $\det(\text{\rm Hess} \lambda(p)) \leq 0.$
It follows that $p_0$ is never a cuspidal lips of $\bx$.
\par
$(2)$ By the proof of $(1)$, $\det(\text{\rm Hess} \lambda(p_0)) <0$ if and only if $\alpha_u(u_0)\neq0$ and $$
(\wm(v_0)t_{12}(p_0)+ \wn(v_0)t_{22}(p_0))\sin\theta(p_0)-(\wm(v_0)t_{11}(p_0) + \wn(v_0)t_{21}(p_0))\cos\theta(p_0)\neq0.
$$
Since ${d\bx}_{p_0}=\walpha(v_0)t_{33}(p_0)\mu(u_0)dv$, we have $\eta=\partial/\partial u$. It follows that $\eta\eta\lambda(p_0)=\lambda_{uu}(p_0)$. Therefore, $\eta\eta\lambda(p_0)\neq0$ if and only if $\alpha_u(u_0)\neq0$ and $-n(u_0)\sin\theta(p_0)+m(u_0)\cos\theta(p_0)\neq0$. Moreover, by Proposition \ref{translation-framed-surface-basic-invariants} and under the same assumptions as in Lemma \ref{lemma-framed-surface}, we have $K^{F}(p_0)=0$ and $H^{F}(p_0)=-(1/2)\walpha(v_0)t_{33}(p_0)(\theta_u(p_0)-\ell(u_0))$. Therefore, $\bx$ to be a front near $p_0$ if and only if $\theta_u(p_0)-\ell(u_0)\neq0$.
\enD


\subsection{$\rm(III)$ The case of $\alpha(u_0) \neq 0, \walpha(v_0) = 0, \mu(u_0) \times \wmu(v_0) = 0$}\label{case3}

By the same calculation of \ref{case2}, $p_0$ is a degenerate singular point of $\bx$ and rank ${d\bx}=1$ at $p_0$ under the condition $\rm(III)$. 
We can prove by the similar calculations of proving of Theorem \ref{criteria-case-II}.

\begin{theorem}\label{criteria-case-III}
Under the assumptions $\rm (III)$, we have the following.
\par
$(1)$ $p_0$ is never a cuspidal lips of $\bx$. 
\par
$(2)$ $p_0$ is a cuspidal beaks of $\bx$ if and only if 
$$
\theta_v(p_0) \neq 0, \ \walpha_v(v_0) \neq 0, \ -n(u_0)\sin \theta(p_0) + m(u_0) \cos \theta(p_0) \neq 0
$$
and $\left(\wm(v_0)t_{12}(p_0)+ \wn(v_0)t_{22}(p_0)\right)\sin\theta(p_0) -
\left(\wm(v_0)t_{11}(p_0) + \wn(v_0)t_{21}(p_0)\right)\cos\theta(p_0) \not=0$. 
\end{theorem}


\begin{remark}{\rm
If $\alpha(u_0) =\walpha(v_0) = 0$ and $\mu(u_0) \times \wmu(v_0) \not= 0$, then 
the translation framed base surface $\bx$ never have cuspidal beaks singular points \cite{Fukunaga-Takahashi-2022}.
}
\end{remark}


\subsection{$\rm(IV)$ The case of $\alpha(u_0) = 0, \walpha(v_0) = 0, \mu(u_0) \times \wmu(v_0) = 0$}\label{case4}

In this case, ${\rm rank}\ d\bx=0$ at $p_0$ under the condition $\rm (IV)$.

\begin{theorem}\label{criteria-case-IV}
Under the assumptions $\rm (IV)$, $p_0$ is never a $D_4^\pm$-singular point of $\bx$. 
\end{theorem}
\demo
By a direction calculation, we have $\lambda_{uu}(p_0)=0, \lambda_{uv}(p_0)=0$ and $\lambda_{vv}(p_0)=0$. Therefore, we have $\text{\rm Hess} \lambda(p_0)=0$. 
Hence, $p_0$ is never a $D_4^\pm$-singular point of $\bx$. 
\enD

\begin{remark}{\rm
If $\alpha(u_0) =\walpha(v_0) = 0$ and $\mu(u_0) \times \wmu(v_0) \not= 0$, then 
the translation framed base surface $\bx$ may have $D^+_4$-singular points \cite{Fukunaga-Takahashi-2022}.
}
\end{remark}

\appendix

\section{Framed curves and regular space curves}\label{framed-curves}

A framed curve in the $3$-dimensional Euclidean space is a smooth space curve with a moving frame, in detail see \cite{Honda-Takahashi-2016}. 
Let $I$ be an interval of $\R$. 
\begin{definition}\label{framed.curve}{\rm
We say that $(\gamma,\nu_1,\nu_2):I \rightarrow \mathbb{R}^3 \times \Delta$ is a {\it framed curve} if $\dot{\gamma}(t) \cdot \nu_1(t)=0$ and $\dot{\gamma}(t) \cdot \nu_2(t)=0$ for all $t \in I$. 
We say that $\gamma:I \to \R^3$ is a {\it framed base curve} if there exists $(\nu_1,\nu_2):I \to \Delta$ such that $(\gamma,\nu_1,\nu_2)$ is a framed curve. 
}
\end{definition}

We denote $\mu(t) = \nu_1(t) \times \nu_2(t)$. 
Then $\{ \nu_1(t),\nu_2(t),\mu(t) \}$ is a moving frame along the framed base curve $\gamma(t)$ in $\R^3$ and we have the Frenet type formula,
$$
\left(
\begin{array}{c}
\dot{\nu_1}(t)\\
\dot{\nu_2}(t)\\
\dot{\mu}(t)
\end{array} \right)=
\left(
\begin{array}{ccc}
0 & \ell(t) & m(t)\\
-\ell(t) & 0 & n(t)\\
-m(t) & -n(t) & 0
\end{array}\right)
\left(
\begin{array}{c}
\nu_1(t)\\
\nu_2(t)\\
\mu(t)
\end{array}\right), 
\ \dot{\gamma}(t)=\alpha(t)\mu(t),
$$
where $\ell(t) = \dot{\nu_1}(t) \cdot \nu_2(t)$, $m(t) = \dot{\nu_1}(t) \cdot \mu(t), n(t) = \dot{\nu_2}(t) \cdot \mu(t)$ and $\alpha(t)=\dot{\gamma}(t) \cdot \mu(t)$. 
We call the mapping $(\ell,m,n,\alpha)$ {\it the (framed) curvature of the framed curve} $(\gamma,\nu_1,\nu_2)$. 
Note that $t_0$ is a singular point of $\gamma$ if and only if $\alpha(t_0) = 0$. 

Suppose that $\gamma:I \to \R^3$ is a regular space curve, that is, $\dot{\gamma}(t) \not=0$ for all $t \in I$, where $\dot{\gamma}(t)=(d\gamma/dt)(t)$. 
We say that $\gamma$ is {\it non-degenerate}, or $\gamma$ satisfies the {\it non-degenerate condition} if $\dot{\gamma}(t) \times \ddot{\gamma}(t) \not=0$ for all $t \in I$. 

If we take the arc-length parameter $s$, that is, $|\gamma'(s)|=1$ for all $s$, then the tangent vector, the principal normal vector and the bi-normal vector are given by
$$
\bt(s)=\gamma'(s), \ \bn(s)=\frac{\gamma''(s)}{|\gamma''(s)|}, \ \bb(s)=\bt(s) \times \bn(s),
$$
where $\gamma'(s)=(d\gamma/ds)(s)$. 
Then $\{\bt(s),\bn(s),\bb(s)\}$ is a moving frame of $\gamma(s)$ and we have the Frenet-Serret formula: 
$$
\left(
\begin{array}{c}
\bt'(s)\\
\bn'(s)\\
\bb'(s)
\end{array}
\right)
=
\left(
\begin{array}{ccc}
0&\kappa(s)&0\\
-\kappa(s)&0&\tau(s)\\
0&-\tau(s)&0
\end{array}
\right)
\left(
\begin{array}{c}
\bt(s)\\
\bn(s)\\
\bb(s)
\end{array}
\right),
$$
where 
$$
\kappa(s)=|\gamma''(s)|, \ \tau(s)=\frac{{\rm det}(\gamma'(s),\gamma''(s),\gamma'''(s))}{\kappa^2(s)}.
$$
The relation between framed curve and non-degenerate curve is as follows. 
Suppose that $\gamma$ is a non-degenerate curve with arc-length parameter and 
the curvature $\kappa$ and the torsion $\tau$. 
Then $(\gamma,\bn,\bb)$ is a framed curve with curvature $(\ell,m,n,\alpha)=(\tau,-\kappa,0,1).$

\section{Generalised framed surfaces and framed surfaces}\label{GFS-FS}

We review the theory of generalised framed surfaces in the Euclidean 3-space, in detail see \cite{Takahashi-Yu}.
Let $(\bx,N_1,N_2):U \to \R^3 \times \Delta$ be a smooth mapping and $U$ be a simply connected domain in $\R^2$.
We denote $\nu=\bx_u \times \bx_v$.

\begin{definition}\label{generalised.framed.surface}{\rm
We say that $(\bx,N_1,N_2):U \to \R^3 \times \Delta$ is a {\it generalised framed surface} if there exist smooth functions $A, B:U \to \R$ such that $\nu(u,v)=A(u,v) N_1(u,v)+B(u,v) N_2(u,v)$ for all $(u,v) \in U$.
We say that $\bx:U \to \R^3$ is a {\it generalised framed base surface} if there exists $(N_1,N_2):U \to \Delta$ such that $(\bx,N_1,N_2)$ is a generalised framed surface.
}
\end{definition}

We denote $N_3(u,v)=N_1(u,v) \times N_2(u,v)$.
Then $\{N_1(u,v),N_2(u,v),N_3(u,v)\}$ is a moving frame along $\bx(u,v)$
and we have the following systems of differential equations:
$$
\begin{pmatrix}
\bx_u \\
\bx_v
\end{pmatrix}
=
\begin{pmatrix}
a_1 & b_1 & c_1\\
a_2 & b_2 & c_2
\end{pmatrix}
\begin{pmatrix}
N_1 \\
N_2 \\
N_3
\end{pmatrix},
$$
$$
\begin{pmatrix}
N_{1u} \\
N_{2u} \\
N_{3u}
\end{pmatrix}
=
\begin{pmatrix}
0 & e_1 & f_1 \\
-e_1 & 0 & g_1 \\
-f_1 & -g_1 & 0
\end{pmatrix}
\begin{pmatrix}
N_1 \\
N_2 \\
N_3
\end{pmatrix}
, \
\begin{pmatrix}
N_{1v} \\
N_{2v} \\
N_{3v}
\end{pmatrix}
=
\begin{pmatrix}
0 & e_2 & f_2 \\
-e_2 & 0 & g_2 \\
-f_2 & -g_2 & 0
\end{pmatrix}
\begin{pmatrix}
N_1 \\
N_2 \\
N_3
\end{pmatrix},
$$
where $a_i,b_i,c_i,e_i,f_i,g_i:U \to \R, i=1,2$ are smooth functions with $a_1b_2-a_2b_1=0$.
We call the functions {\it basic invariants} of the generalised framed surface.
We denote the above matrices by $\mathcal{G}, \mathcal{F}_1, \mathcal{F}_2$, respectively.
We also call the matrices $(\mathcal{G}, \mathcal{F}_1, \mathcal{F}_2)$ {\it basic invariants} of the generalised framed surface $(\bx,N_1,N_2)$.
By definition, we have
$$
A(u,v)={\det}
\begin{pmatrix}
b_1(u,v) & c_1(u,v)\\
b_2(u,v) & c_2(u,v)
\end{pmatrix}, \
B(u,v)=-{\det}
\begin{pmatrix}
a_1(u,v) & c_1(u,v)\\
a_2(u,v) & c_2(u,v)
\end{pmatrix}.
$$

Since the integrability conditions $\bx_{uv}=\bx_{vu}$ and $\mathcal{F}_{2u}-\mathcal{F}_{1v}
=\mathcal{F}_1\mathcal{F}_2-\mathcal{F}_2\mathcal{F}_1$, 
the basic invariants should be satisfied some conditions.

Note that there are fundamental theorems for generalised framed surfaces, namely,
the existence and uniqueness theorems for the basic invariants of generalised framed surfaces (cf. \cite{Takahashi-Yu}).

We also review the theory of framed surfaces in the Euclidean 3-space, in detail see \cite{Fukunaga-Takahashi-2019, Fukunaga-Takahashi-2020}.
Let $(\bx,\bn,\bs):U \to \R^3 \times \Delta$ be a smooth mapping.
\begin{definition}\label{framed.surface}{\rm
We say that $(\bx,\bn,\bs):U \to \R^3 \times \Delta$ is a {\it framed surface} if $\bx_u (u,v) \cdot \bn (u,v)=\bx_v(u,v) \cdot \bn(u,v)=0$ for all $(u,v) \in U$, where $\bx_u(u,v)=(\partial \bx/\partial u)(u,v)$ and $\bx_v(u,v)=(\partial \bx/\partial v)(u,v)$.
We say that $\bx:U \to \R^3$ is a {\it framed base surface} if there exists $(\bn,\bs):U \to \Delta$ such that $(\bx,\bn,\bs)$ is a framed surface.
}
\end{definition}
By definition, the framed base surface is a frontal.
The definition and properties of frontals see \cite{Arnold1,Arnold2}.
On the other hand, the frontal is a framed base surface at least locally.

We denote $\bt(u,v)=\bn(u,v) \times \bs(u,v)$.
Then $\{\bn(u,v),\bs(u,v),\bt(u,v)\}$ is a moving frame along $\bx(u,v)$
and we have the following systems of differential equations:
$$
\begin{pmatrix}
\bx_u \\
\bx_v
\end{pmatrix}
=
\begin{pmatrix}
a_1 & b_1 \\
a_2 & b_2
\end{pmatrix}
\begin{pmatrix}
\bs \\
\bt
\end{pmatrix},
$$
$$
\begin{pmatrix}
\bn_u \\
\bs_u \\
\bt_u
\end{pmatrix}
=
\begin{pmatrix}
0 & e_1 & f_1 \\
-e_1 & 0 & g_1 \\
-f_1 & -g_1 & 0
\end{pmatrix}
\begin{pmatrix}
\bn \\
\bs \\
\bt
\end{pmatrix}
, \
\begin{pmatrix}
\bn_v \\
\bs_v \\
\bt_v
\end{pmatrix}
=
\begin{pmatrix}
0 & e_2 & f_2 \\
-e_2 & 0 & g_2 \\
-f_2 & -g_2 & 0
\end{pmatrix}
\begin{pmatrix}
\bn \\
\bs \\
\bt
\end{pmatrix},
$$
where $a_i,b_i,e_i,f_i,g_i:U \to \R, i=1,2$ are smooth functions.
We call these functions {\it basic invariants} of the framed surface.
We denote the above matrices by $\mathcal{G}, \mathcal{F}_1, \mathcal{F}_2$, respectively.
We also call the matrices $(\mathcal{G}, \mathcal{F}_1, \mathcal{F}_2)$  {\it basic invariants} of the framed surface $(\bx,\bn,\bs)$.
Note that $(u,v)$ is a singular point of $\bx$ if and only if ${\rm det}\ \mathcal{G}(u,v)=0$.
We define the curvature $C^F=(J^F,K^F,H^F)$ of the framed surface $(\bx,\bn,\bs)$ by
\begin{align*}
&J^F={\det}
\begin{pmatrix}
a_1 & b_1\\
a_2 & b_2
\end{pmatrix}, \ 
K^F={\det}
\begin{pmatrix}
e_1 & f_1\\
e_2 & f_2
\end{pmatrix}, \ 
H^F=-\frac{1}{2} 
\left\{
\det \begin{pmatrix}
a_1 & f_1\\
a_2 & f_2
\end{pmatrix}
- \det \begin{pmatrix}
b_1 & e_1\\
b_2 & e_2
\end{pmatrix}
\right\}.
\end{align*}
Since the integrability conditions $\bx_{uv}=\bx_{vu}$ and $\mathcal{F}_{2,u}-\mathcal{F}_{1,v}
=\mathcal{F}_1\mathcal{F}_2-\mathcal{F}_2\mathcal{F}_1$,
the basic invariants should be satisfied some conditions.
Note that there are fundamental theorems for framed surfaces, namely,
the existence and uniqueness theorems for the basic invariants of framed surfaces (cf. \cite{Fukunaga-Takahashi-2019}).
We give a relation between generalised framed base surfaces and framed base surfaces.

\begin{theorem}[\cite{Takahashi-Yu}]\label{FBS.condition}
Let $(\bx,N_1,N_2):U \to \R^3 \times \Delta$ be a generalised framed surface with $\nu=A N_1+B N_2$.
\par
$(1)$ If $\bx$ is a framed base surface, then the functions $A$ and $B$ are linearly  dependent.
\par
$(2)$ Suppose that the set of regular points of $\bx$ is dense in $U$.
If the functions $A$ and $B$ are linearly dependent, then $\bx$ is a framed base surface at least locally.
\end{theorem}

\begin{proposition}\label{Legendre.immersion}
Let $(\bx,\bn,\bs) : U \to \R^3 \times \Delta$ be a framed surface and $p \in U$. 
\par
$(1)$ Suppose that ${\rm rank} (d\bx)=1$ at $p$. 
Then $(\bx,\bn): U \to \R^3 \times S^2$ is a Legendre immersion around $p$ if and only if $H^F(p) \neq 0$.
\par
$(2)$ Suppose that ${\rm rank} (d\bx)=0$ at $p$. 
Then $(\bx,\bn): U \to \R^3 \times S^2$ is a Legendre immersion around $p$ if and only if $K^F(p) \neq 0$.
\end{proposition}

\section{Criterion for singularities}\label{criterion}

We introduce singular points and their criterion which are used in this paper. 
\bigskip

If a map-germ $f : (\R^2, 0) \to (\R^3, 0)$ satisfies rank $df_0 = 1$, the singular point $0$ is called
corank one. 
If $f : (\R^2, 0) \to (\R^3, 0)$ has a corank one singular point at $0$, then there exist vector fields $(\xi, \eta)$ near the origin such that $df_0(\eta_0) = 0$ and $\xi_0, \eta_0 \in T_0 \R^2$ are linearly
independent. 
We define a function $\varphi : (\R^2, 0) \to \R$ by $\varphi= {\rm det}(\xi f, \eta f, \eta \eta f),$ 
where $\zeta g : (\R^2, 0) \to (\R^3, 0)$ is the directional derivative of a vector valued function $g$ by a vector field $\zeta$. 

We say that a singular point $p$ of a mapping $\bx : U \rightarrow \mathbb{R}^3$ is a
{\it cross cap} $S_0$ (respectively, Chen Matsumoto Mond $+$, Chen Matsumoto Mond $-$ singular point $S_{1}^{\pm}$)
if $\bx$ is $\mathcal{A}$-equivalent, that is, equivalent by diffeomorphisms of 
the source and of the target to the map germ $(u,v) \mapsto (u,uv,v^2)$ (respectively, $(u,v) \mapsto (u, v^2, v(u^2 \pm v^2)))$) at $p$ (see \cite{Mond}).
In \cite{Saji2010}, K.Saji gave the criteria for $S_{1}^{\pm}$ singular points.

\begin{theorem}[\cite{Saji2010}]\label{CMM}
Let $\bx : (\mathbb{R}^2,0) \rightarrow (\mathbb{R}^3,0)$ be a map-germ and $0$ be a corank $1$ singular point. 
\par
$(1)$ $\bx$ at $0$ is a cross cap $S_0$ singular point if and only if $\xi \varphi \not= 0$ at $0$, that is, $d \varphi \not= 0$ at $0$.
\par
$(2)$ $\bx$ at $0$ is a $S_{1}^+$ singular point if and only if $\varphi$ has a critical point at $0$, $\det(\text{\rm Hess} \varphi(0))<0$ and two vectors $\xi\bx(0)$ and
 $\eta\eta\bx(0)$ are linearly independent.
\par
$(3)$ $\bx$ at $0$ is a $S_{1}^-$ singular point if and only if $\varphi$ has a critical point at $0$ and $\det(\text{\rm Hess} \varphi(0))>0$. 
Here,  $\text{\rm Hess} \varphi(0)$ means the Hessian matrix of $\varphi$ at $0$.
\end{theorem}

We introduce cuspidal edges, swallowtails and cuspidal cross caps
which are known as the generic singular points of a frontal from
$\mathbb{R}^2$ to $\mathbb{R}^3$. 
We say that a singular point $p$ of a mapping $\bx : U \rightarrow \mathbb{R}^3$ is a
{\it cuspidal edge} (respectively, {\it swallowtail} or {\it cuspidal cross cap})
if $\bx$ is $\mathcal{A}$-equivalent to the map germ $(u,v) \mapsto
(u,v^2,v^3)$ (respectively, $(u,v) \mapsto (u, 4v^3+2uv, 3v^4+uv^2)$ or $(u,v) \mapsto (u, v^2, uv^3)$) at $p$.

We recall the criteria for singular points of frontals stated in
\cite{Fujimori-Saji-Umehara-Yamada, Kokubu-Rossman-Saji-Umehara-Yamada} (see also, \cite{Izumiya-Saji}).
Let $\bx : U \rightarrow \mathbb{R}^3$ be the frontal of a Legendre
surface $(\bx,\bn)$. 
We define a function $\lambda : U \rightarrow
\mathbb{R}$ by $\lambda(u,v) = \det(\bx_u,\bx_v,\bn)(u,v)$ where $(u,v)$
is a coordinate system on $U$. 
We call the function $\lambda$ a {\it discriminant function} (or, a {\it signed area density function}).
When a singular point $p$ of $\bx$ is non-degenerate, that is, $d\lambda(p) \neq 0$, there exists a smooth parametrisation $\delta : (-\varepsilon,\varepsilon) \rightarrow U$, $\delta(0)=p$ of the singular set $S(\bx)$. 
We call the curve $\delta$ the singular curve of $\bx$. 
Moreover, there exists a smooth vector field $\eta(t)$ along $\delta(t)$ satisfying that $\eta(t)$ generates $\ker d\bx_{\delta(t)}$. 
Now we define a function $\phi_{\bx} : (-\varepsilon,\varepsilon) \rightarrow \mathbb{R}$
by $\phi_{\bx}(t)=\det((\bx \circ \delta)^{\prime}, \bn \circ \delta,
d\widehat{\bn}(\eta))(t)$ where $\widehat{\bn}$ is a normal vector
field of $\bx$ which does not have zero points. 
By using these notations, we have the following theorem.

\begin{theorem}[\cite{Fujimori-Saji-Umehara-Yamada,Kokubu-Rossman-Saji-Umehara-Yamada}] \label{criteria_Izumiya-Saji}
Let $(\bx,\bn) : U \rightarrow \mathbb{R}^3 \times S^2$ be a Legendre surface and $p \in U$ be a non-degenerate singular point of $\bx$.
Then the following assertions hold.
\par
$(1)$ If $\eta\lambda(p) \neq 0$, then $\bx$ to be a front near $p$ if
 and only if $\phi_{\bx}(0) \neq 0$ holds.
\par
$(2)$ The map germ $\bx$ at $p$ is a cuspidal edge if and only if $\bx$ to be a front near $p$ and
 $\eta\lambda(p) \neq 0$ hold.
\par
$(3)$ The map germ $\bx$ at $p$ is a swallowtail if and only if $\bx$ to be a front near $p$ and
 $\eta\lambda(p) = 0$ and $\eta\eta\lambda(p) \neq 0$ hold.
\par
$(4)$ The map germ $\bx$ at $p$ is a cuspidal cross cap if and only if $\eta \lambda(p) \neq 0$,
 $\phi_{\bx}(0) = 0$ and $\phi_{\bx}^{\prime}(0) \neq 0$ hold.
\par
Here, $\eta\lambda : U \rightarrow \mathbb{R}$ means the directional derivative of $\lambda$ by the vector field $\tilde\eta$, where $\tilde\eta$ is an extended vector field of $\eta$ to $U$. 
\end{theorem}

We consider cuspidal lips, cuspidal beaks and $D_4^{\pm}$-singular points.
They are known as the generic singular points of one-parameter families of fronts
(see \cite{Arnold1}). 
Let $\bx : (U,p) \rightarrow \mathbb{R}^3$ be a map germ. 
We call $p$ is a cuspidal lips (respectively, cuspidal beaks) if $\bx$ is 
$\mathcal{A}$-equivalent to the map germ $(u,v) \mapsto
(3u^4+2u^2v^2,u^3+uv^2,v)$ 
(respectively, $(u,v) \mapsto (3u^4-2u^2v^2,u^3-uv^2,v)$). 
Moreover, we call $p$ is a $D_4^{+}$ singular point (respectively,
a $D_4^{-}$ singular point) if
$\bx$ is $\mathcal{A}$-equivalent to the map germ
$(u,v) \mapsto (uv, u^2 + 3 v^2, u^2 v + v^3)$
(respectively, $(u,v) \mapsto (uv, u^2 - 3 v^2, u^2 v - v^3)$) at $p$.
Criteria for cuspidal lips and cuspidal beaks are first given by
S.Izumiya, K.Saji and the third author in \cite{IST2010}. 
In this paper, we follow the notations in \cite{Fukui-Hasegawa2012}.

\begin{theorem}[\cite{Fukui-Hasegawa2012}]\label{CLCB}
Let $\bx : U \rightarrow \mathbb{R}^3$ be a front and $p \in U$ be a
 degenerate singular point of $\bx$.

$(1)$ The germ of the front $\bx$ at $p$ is a cuspidal lips if and only if ${\rm rank}(d\bx(p)) = 1$ and
 $\det(\text{\rm Hess} \lambda(p))>0$.

$(2)$ The germ of the front $\bx$ at $p$ is a cuspidal beaks if and only if ${\rm rank}(d\bx(p)) = 1$,
 $\det(\text{\rm Hess} \lambda(p))<0$ and $\eta\eta\lambda(p) \neq 0$.
\end{theorem}

The criteria for $D_4^{\pm}$ singular points are given by K.Saji in \cite{Saji2}.
\begin{theorem}[\cite{Saji2}]\label{D_4_criteria}
Let $\bx : (\mathbb{R}^2,0) \rightarrow (\mathbb{R}^3,0)$ be a front 
and $(\bx, \bn)$ its Legendrean lift. The germ $\bx$ at $0$ is a $D_4^{+}$
 singular point (respectively, $D_4^{-}$ singular point) if and only if the following two conditions hold.
\par
 $(a)$ The rank of the differential map $d\bx(0)$ is equal to zero, 
\par
$(b)$ $\det(\text{\rm Hess} \lambda(0))<0$ (respectively, $\det(\text{\rm Hess}
 \lambda(0))>0$).
\end{theorem}


Tomonori Fukunaga, 
\par
Fukuoka Institute of Technology, Fukuoka 811-0295, Japan,
\par
E-mail address: fukunaga@fit.ac.jp

\bigskip
Nozomi Nakatsuyama, 
\par
Muroran Institute of Technology, Muroran 050-8585, Japan,
\par
E-mail address: 25096009b@muroran-it.ac.jp

\bigskip
Masatomo Takahashi, 
\par
Muroran Institute of Technology, Muroran 050-8585, Japan,
\par
E-mail address: masatomo@muroran-it.ac.jp
\end{document}